 \documentclass[final,1p,times]{elsarticle}

\usepackage{lscape}

\usepackage{amssymb}
 \usepackage{amsthm}

\newcommand\GL{\operatorname{GL}}
\newcommand\Sp{\operatorname{Sp}}
\newcommand\Orth{\operatorname{O}}
\newcommand\PSL{\operatorname{PSL}}
\newcommand\U{\operatorname{U}}

\usepackage{calrsfs}

\usepackage{hyperref}
\usepackage{mathtools}

\makeatletter
\newcounter{@myftn}
\def\myfootnote{\@ifnextchar [\@xfootnote {\stepcounter {@myftn}\protected@xdef \@thefnmark%
{\alph{@myftn}}\@footnotemark \@footnotetext}}
\makeatother

\journal{Indagationes Math}

\begin{document}

\begin{frontmatter}



\title{On symplectic transformations}


\author{T.A. Springer}
\ead{W.vanderKallen@uu.nl [editor]}
\affiliation{organization={Universiteit Utrecht},
            addressline={Postbus 80.010}, 
            city={Utrecht},
            postcode={3508 TA}, 
            country={Nederland}}

\begin{abstract}
This is an English translation\myfootnote{\emph{Prepared and edited by Wilberd van der Kallen. Translated and published by permission of the heirs of T.A.~Springer. 
The editor has corrected some minor glitches in translating from the Dutch original; 
these are indicated by footnotes in italic, with an alphabetical label. 
The numbered footnotes are due to the author of the thesis.}} of the PhD thesis `Over symplectische transformaties' 
that Tonny Albert Springer, `born in 's-Gravenhage in
1926',  submitted as thesis for---as is stated on the original 
frontispiece---\emph{the degree of doctor in mathematics and physics at Leiden 
University on the authority of the rector magnificus Dr.~J.~H.~Boeke, professor in 
the faculty of law, to be defended against the objections of the
Faculty of Mathematics and Physics  on Wednesday October 17
1951 at 4~p.m.}, with promotor Prof.\ dr.\ H.~D.~Kloosterman.
\end{abstract}




\end{frontmatter}




%


\section*{Reproduction of the original frontispiece of the thesis} 

\begin{center}
\includegraphics[width=10cm]{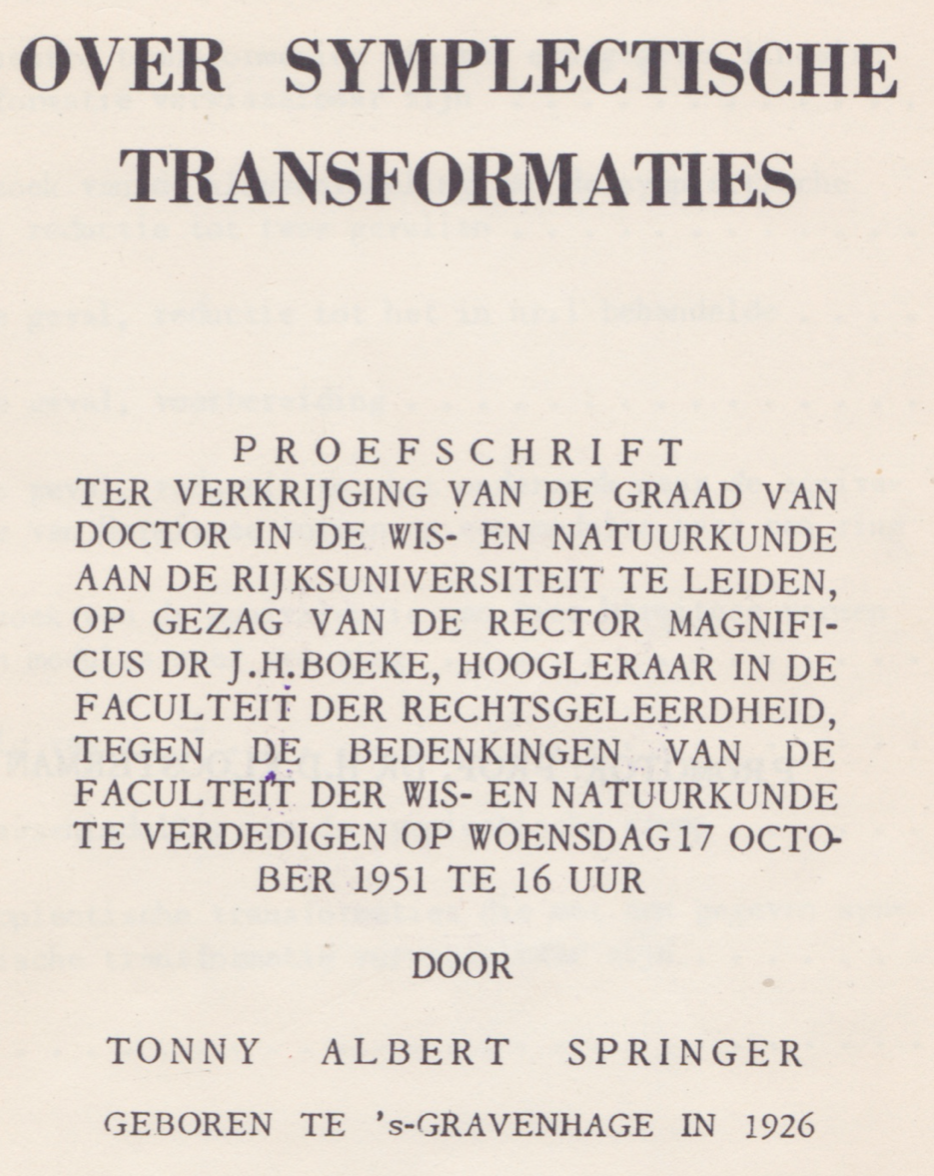}
\end{center}

\section*{Introduction}

The classification of conjugacy classes of  elements of the general
linear group $\GL_n(K)$, the group of invertible linear
transformations in $n$ variables over a commutative coefficient field K, is easily found with the help of the
results from the theory of canonical forms of linear
transformations\footnote{For this
theory, see, e.g.\ \cite[ p.~121]{vdW}.}. In this thesis we investigate the
classification of conjugacy classes in the symplectic group $\Sp_n(K)$, which is the subgroup
of $\GL_n(K)$ consisting of all linear transformations that leave invariant a given
skew-symmetric form. We assume that the
characteristic of $K$ differs from $2$. 

It is found that, as with the
group $\GL_n(K)$, every conjugacy class of $\Sp_n(K)$ can be characterized by
certain invariants. These are irreducible polynomials and systems of
integers (such invariants also occur for $\GL_n(K)$), and, 
additionally, equivalence classes of certain Hermitian or quadratic
forms (see section~\ref 9). 

For the sake of completeness,  the classification of conjugacy classes in
$\GL_n(K)$ is
treated first. Also, the structure of the normalizer of
an element of $\Sp_n(K)$ is investigated.

It seems that the classification of conjugacy classes in $\Sp_n(K)$
for an arbitrary field $K$ has not been treated before.\myfootnote{\emph{i.e., before 1951.}} Special
cases have been discussed by L.~E.~Dickson \cite{LED1}, \cite{LED2} ($\Sp_4(K)$ and $\Sp_6(K)$ if $K$ is a
finite field) and by J.~Williamson \cite{W} (if $K$ is the field of real
numbers). 

The classification of conjugacy classes in the other classical linear
groups (the orthogonal and unitary groups) can be examined by an analogous method. We hope to return to that.   

\section{The classification of conjugacy classes in the
general linear group.}\label{1}
 We assume to be given a
vector space $E$ of finite dimension $n$ over the commutative field $K$.
Suppose $u$ is a linear transformation of $E$. To each element $f =
\sum_{\nu=0}^kY_\nu X^\nu$ ($Y_\nu \in K$) of the 
ring $K[X]$ of polynomials in one variable $X$ over $K$ we assign
the linear transformation $f(u) =\sum_{\nu=0}^kY_\nu u^\nu$ of $E$. This assignment
is a homomorphism from $K[X]$ onto the ring of linear transformations $f(u)$.
If $d$ is the unique polynomial of lowest possible degree
and with leading coefficient $1$ for which $d(u) = 0$ holds ($d$ is called the
minimal polynomial of $u$), then the ring consisting of the $f(u)$ is isomorphic to  the
residue class ring $R=K[X]/(d)$.

 Now in $K [X]$, $d$ can be written  as
product of mutually distinct irreducible polynomials with highest
coefficient $1$.
Suppose this decomposition is $d = \prod_{i=1}^s p_i^{k_i}$.  Denote by  $E_i$ the subspace of 
$E$, which consists of all $x\in E$ for which $p_i^{k_i}(u)(x)=0$ holds. Every $E_i$
is transformed into itself by $u$ (because $up_i^{k_i}(u)=p_i^{k_i}(u)u$). 
Each element of $E$ can be written as a sum of elements 
from $E_i$ (there are polynomials $g_i$ such that $\sum_{i=1}^s g_i\frac{d}{p_i^{k_i}}=1$ and for an
$x\in E$ one thus has $x=\sum_{i=1}^sx_i$ with $x_i=g_i\frac{d}{p_i^{k_i}}(u)(x)\in E_i$).
This can be done in one way only:
if $\sum_{i=1}^sx_i=0$ with $x_i\in E_i$, then $\frac{d}{p_i^{k_i}}(u)(x_i)=0$
and, as $\frac{d}{p_i^{k_i}}$ and $p_i^{k_i}$ are relatively prime, there exist polynomials $a$ and $b$ such 
that $a\frac{d}{p_i^{k_i}}+ bp_i^{k_i} =1$, from which it follows that $x_i= (a\frac{d}{p_i^{k_i}}+ bp_i^{k_i})(u)(x_i) =0$. We 
see that $E$ is the direct sum of the $E_i$ ($1\leq i\leq s$). The restriction of $u$ to 
$E_i$ is a linear transformation $u_i$ of $E_i$ with minimal polynomial $p_i^{k_i}$. 
For $x =\sum_{i=1}^sx_i$ with $x_i\in E_i$ one has $u(x)=\sum_{i=1}^su_i(x_i)$. We investigate the individual $u_i$ ($1\leq i\leq s$).

We can assume that $u$ itself is a linear
transformation with minimal polynomial $p^k$, where $p$ is an irreducible
polynomial of degree $g$.\footnote{When we speak about an
irreducible polynomial, we shall here and in the sequel assume that the leading 
coefficient of this polynomial is $1$.} 
Since each element of $R = K[X]/(p^k)$ yields a
transformation $f(u)$, we can give $E$ the structure of an $R$-module.
We call this module $M$. The multiplication with elements of
$R$ is thus defined as follows: if $\kappa$ is the canonical homomorphism from $K[X]$ onto
$R$, then for $\rho= \kappa(f)\in R$ one has: $\rho x =
 f(u)(x)$.

The ring $R$ is commutative with unit element and has a simple 
structure: when $\pi= \kappa(p)$, every nonzero element of $R$ can be written
as $\pi^a\epsilon$ \ $(0\leq a\leq k)$, where $\epsilon$ is a unit (invertible element) in $R$.\ \ 
 The radical of $R$ is $(\pi)$. Also, every element $\rho$
of $R$ can be unambiguously written as
$\rho=\sum_{\nu=0}^{kg-1} \alpha_\nu\xi^\nu$ \ $(\alpha_\nu\in K$), where $\xi=\kappa(X)$.
We denote by $R_i$ the residue class ring $R/(\pi^i)$ \ $(0\leq i\leq k)$. 
We
can then consider $R_i$ as a ring with operators from $R$. Also, as a ring, $ R_i$ is isomorphic 
to $K[X]/(p^i)$, in particular, $R_1$ (the residue class ring 
with respect to the radical) is isomorphic to the field $L =K[X]/(p)$. 

We now investigate
the $R$-module $M$ in more detail. We will prove:

\begin{quote}
\emph{There are elements $e_i^j$ \ $( 1\leq i\leq k, \ 1\leq j\leq a_i )$ of $M$ such that $M$ is the direct sum 
of the modules $Re_i^j$.} 
\end{quote}

(As usual, by $Re_i^j$ is meant the
submodule of $M$ consisting of the elements $\rho e_i^j$ with $\rho\in R$). 

We define the $e_i^j$ as follows. Take as $e_k^j$ \ $(1\leq j\leq a_k)$ a maximal 
system of elements of $M$ such that $\pi^{k-1}e_k^j\neq0$ and such that every module 
$Re_k^j$ has only the zero element in common with the sum of the others. Suppose that
the $e_i^j$ \ $(n<i\leq k)$ are already known. If $n>0$, take as $e_n^j$ \ $( 1\leq j\leq a_n)$ a maximal
system of elements of $M$ such that $\pi^n e_n^j=0$, $\pi^{n-1}e_n^j\neq0$, and such that
each module  $Re_n^j$  \ $(n\leq i\leq k)$ has only the zero element
in common with the sum of the other $Re_i^j$ \ $(n\leq i\leq k)$. 
 It is of course possible that such $e_n^j$ do not exist. We then
set $a_n= 0$ (however, $a_k\neq0$).  We now prove that every $x\in M$ is a sum
of elements from the $Re_i^j$.  Take an $x\in M$. Suppose $\pi^hx = 0$ and $\pi^{h-1}x\neq0$ if 
$h>0$.\myfootnote{\emph{edited to make the case $h=0$ meaningful.}}
The assertion is proved by complete induction on $h$. For $h=
0$ it is correct. Assume, that the correctness is already proved for all $n<h$.
Due to the definition of the $e_h^j$, there are elements $\rho$ and $\rho_i^j$ of $R$ such
that $\rho x -\sum_{i=h}^k\sum_{j=1}^{a_i}\rho_i^je_i^j=0$
and such that not all $\rho_i^je_i^j$ are zero.  Multiplying 
by a unit of $R$, one can arrange that
$\rho = \pi^n$ \ $(0\leq n<h)$. We assume this now. 
Since $\pi^{ h-n} \rho x = 0$, we have $\sum_{i=h}^k\sum_{j=1}^{a_i}\pi^{h-n}\rho_i^je_i^j=0$, which is possible only if $\pi^{h-n}\rho_i^j$
 is a multiple of $\pi^i$. We can therefore state that $\rho_i^j=\pi^n\sigma_i^j$.
 It follows that $\pi^n(x-\sum_{i=h}^k\sum_{j=1}^{a_i}\sigma_i^je_i^j) = 0$. However, according to the inductive assumption, the element in parentheses
is a sum of elements from the
$Re_i^j$. The same is therefore true for $x$. Thus $M$ is the sum of the
$Re_i^j$. That it
 is a direct sum of the 
$Re_i^j$ follows from the definition of the $e_i^j$.

 We 
further note that the annihilator of $e_i^j$ is the ideal $(\pi^i)$ \ $(1\leq i\leq k)$, which
shows that $Re_i^j$ is  isomorphic to  the $R$-module $R_i$. Denoting $R_i^{a_i}$ the
direct sum of $a_i$ modules that are  isomorphic to   $R_i$, we see that the
result can also be pronounced like this: 

\begin{quote}
\emph{$M$ is the direct sum of modules  isomorphic to  $R_i^{a_i}$
\  $(1\leq i\leq k)$.}
  \end{quote} 

 It follows from the above that the
vector space $E$ over $K$ has a basis consisting of elements $u^h(e_i^j)$ \ $(0\leq h<ig)$, and that 
with respect to this basis $u$ is represented by a matrix
which is completely 
determined by $p$, $K$ and the $a_i$. 
Here $g\sum_{i=1}^kia_i=n$. 
One thus arrives at the well-known canonical forms of matrices\footnote{The derivation given here is 
essentially the one of \cite[p. 307]{SvdW}}. 

If $E_1$
is a second $n$-dimensional vector space over $K$, $u_1$ is a linear
transformation of $E_1$ and $t$ a linear mapping from $E$ onto $E_1$ such that 
$tu(x) = u_1t(x)$ \ $(x\in E)$, then it is easy to see that to $u_1$
belong the same $p$, $K$ and $a_i$ as to $u$. In the special case that $E_1=E$  
and that $t$ is a linear mapping from $E$ onto itself, it follows that
the same 
$p$, $K$ and $a_i$ belong to $u_1=tut^{-1}$ as belong to $u$. The other way around: if the linear transformations $u$ and
$u_1$ of $E$ with minimal polynomial $p^k$ have the same $ a_i$, then there are two
systems of basis vectors of $E$ with respect to which $u$ and $u_1$  are represented by the same
matrix. The transition from one basis to the other
then gives a $t$ for which $u_1 = tut^{-1}$.

Finally, let there be given an irreducible polynomial $p$ of degree $g$,
integers $k>0$ and $a_i \geq0$ \ $( 1\leq i\leq k,\ a_k>0)$ such that $g\sum_{i=1}^kia_i=n$. Then
multiplication by the element $\xi$ of $R$ gives a linear  transformation $u'$
of the direct sum of the modules $R_i^{a_i}$ \ $(1\leq i\leq  k)$. This direct sum is a vector space over $K$. The linear  transformation $u'$, seen as a linear transformation of this
vector space, has as 
minimal polynomial $p^k$ and the numbers $a_i$ belong to it. Since two $n$-dimensional
vector spaces over $K$ can already be mapped onto each other linearly,
it follows that there are linear transformations $u$ of $E$ with given
$a_i$ (for fixed $p$ and $k$). As far as the
classification of conjugacy classes in the group $\GL_n (K)$ of invertible linear
transformations of $E$ is concerned, it follows from the foregoing that: 

\begin{quote}
\emph{Each conjugacy class of $\GL_n(K)$ is unambiguously determined by a number of mutually
distinct irreducible polynomials $p_i$ different from $X$, with degrees
$g_i$, by integers $k_i>0$ and integers $ a^i_j\geq0$ \ $( 1\leq i\leq s,\  1\leq j\leq k_i,\ 
a^i_{k_i} >0)$ such that $\sum_{i=1}^sg_i\sum_{ j= 1}^{k_i} ja^i_j = n$.}
\end{quote} 

\section{The linear transformations that commute with a given one.}\label{2} 
We assume to be given a linear transformation $u$ of the
vector space $E$ with minimal polynomial $d$ and we investigate the linear transformations
of $E$ that commute with $u$. The vector space $ E$ can again be decomposed
into the subspaces $E_i$ \ $(1\leq i\leq s)$. If $v$ is a linear transformation of $E$ 
commuting with $u$, then $p_i^{k_i}(u)(v(x))=v(p_i^{k_i} (u)(x))$ \ $(x\in E)$.
It follows  that $v$ transforms every $E_i$ into itself. We can
therefore again restrict ourselves to the case where $u$ has a minimal polynomial $p^k$.

To a $v$ commuting with $u$ belongs a self-map of the $R$-module $M$,
which we shall also call $v$. When $\rho = \kappa( f) \in R$, we have $v(\rho x) =
v(\kappa(f)x) =v(f(u)(x)) = f(u) (v(x)) = \rho v(x)$. Since $v(x + y) = v(x) +
v(y)$, the self map $v$ of $M$ is an $R$-linear map. Conversely, if
for all $\rho\in R$, we have $v(\rho x) = \rho v(x)$, then also $vu(x) = v(\xi x)=\xi v(x)=uv(x)$
for all $x\in E$. We see: 

\begin{quote}
\emph{To every linear
transformation of the vector space $E$ that commutes with $u$ 
belongs  an $R$-linear mapping
of the $R$-module $M$ into itself, and vice versa.}
\end{quote}

 It follows that the
group invertible linear transformations of $E$ that commute with $u$ is isomorphic
with the group $\GL(M,R)$ of automorphisms of $M$ ($R$-linear maps
of $M$ onto itself). We will examine this group in greater detail.

\

 First we note
 that an endomorphism $v$ of $M$ (i.e., an $R$-linear map of $M$ into
itself) is completely determined by the $v(e_i^j)$. If
\begin{equation}\label{eq1}
 v( e_i^j) = \sum_{p=1}^k\sum_{q=1}^{a_p}\rho_{ip}^{jq}e_p^q\quad(\rho_{ip}^{jq}\in R),
\end{equation}
then $v(\pi^i e_i^j) = 0$ must hold, from which it follows that $\pi^i\rho_{ip}^{jq}\equiv 0\pmod {\pi^p}$
in $R$. Conversely, if one has elements $\rho_{ip}^{jq}$ of $R$, for which these
congruences hold, then  (\ref{eq1})  determines an endomorphism of $M$. Now
$v$ is an automorphism of $M$ if and only if  $v(x) =
0$ implies $x = 0$. This follows directly from the corresponding assertion for a
vector space. 

\

Denote by  $M_i$ \ $(0\leq i\leq k)$ the subspace of $M$ which consists of the
$x\in M$ for which $\pi^i x = 0$. Since $\pi^{k-i}x\in M_i$ for every $x\in M$,
one may view $M/M_i$ as a module  over $R_{k-i}$. Denote by   $\pi^iM$ \ $(0\leq i\leq k)$ the
submodule of $M$ consisting of the elements $\pi^ix$ \ $(x\in M)$. Then $M/\pi^iM$ can
be understood as a module  over $R_i$. Further $ M_i$ is the direct sum of the modules
$Re_p^q$ \ $(1\leq p\leq i)$, \ $R(\pi^{p-i}e_p^q)$ \ $(i+1\leq p\leq k)$. When $v$ is an endomorphism of $M$
one has  $v(\pi^ix) = \pi^iv(x)$ for all $x\in M$, showing that $M_i$ and $\pi^iM$
are transformed by $v$ into themselves. When $v$ is an automorphism, $M_i$ and $\pi^iM$ are mapped onto themselves. 
An endomorphism $v$ of $M$
induces an endomorphism of the $R_{k-i}$-module $M/M_i$. Thus one gets
a homomorphism from $\GL(M,R)$ to  $\GL(M/M_i ,R_{k-i} )$. We will show
that this is a homomorphism onto
$\GL(M/M_i, R_{k-i} )$. 

Since $ M/M_i$ is isomorphic to $\pi^iM$ (isomorphism theorem) and since
$\pi^iM$ is the direct sum of the modules $R(\pi^ie_p^q)$ \ $(i+1\leq p\leq k)$, the $R_{k-i}$-module $M/M_i$ is 
the direct sum of modules $R_{k-i} \bar e_p^q$ \ $( i+1\leq p\leq k,\ 1\leq q\leq a_p)$; if $\psi_i$ is the
 canonical homomorphism from $M$ onto $M/M_i$, then $\bar e_p^q=\psi(e_p^q)$. Let $\phi_{k-i}$ denote
the canonical homomorphism from $R$ onto $R_{k-i}$. 
Suppose that for an
automorphism $w$ of $M/M_i$ it holds that $w(\bar e_p^q) = \sum_{r=i+1}^k\sum_{s=1}^{a_r}\phi_{k-i}(\rho_{pr}^{qs})\bar e_r^s$.
Then  $v(e_p^q) = e_p^q$ \ $(1\leq p\leq i)$, \ $ v(e_p^q) = \sum_{r=i+1}^k\sum_{s=1}^{a_r}\rho_{pr}^{qs}e_r^s$ \  $(i+1\leq p\leq k)$
 determines an endomorphism $v$ of $M$ which on $M/M_i$ induces the
automorphism $w$. Suppose that $v(x) = 0$ for an $x\in M$. Since
$w$ is an 
automorphism, $x$ lies in $M_i$. Write $x=\sum_{p=1}^i\sum_{q=1}^{a_p}\xi_p^qe_p^q+\pi^\alpha y$  \
$(\xi_p^q\in R,\ y\in\sum_{p=i+1}^k\sum_{q=1}^{a_p} Re_p^q)$.
Say $y\in M_n$, $y\notin M_{n-1}$ \ $(n>i)$. From $v(x) = 0$ it follows that $\xi_p^qe_p^q=0$, $\pi^\alpha v(y)=0$.
Since $w$ is an automorphism of $M/M_i$, one also has $v(y)\in M_n$,
$v(y)\notin M_{n-1}$. From $\pi^\alpha v(y)=0$ it then follows that
$\alpha\geq n$, 
so that $\pi^\alpha y=0$.   
Then $x = 0$, so $v$ is an automorphism of $M$.
We have thus proved: 

\begin{quote}
\emph{The assignment that to an automorphism of $M$ associates 
the automorphism it induces on $M/M_i$  gives a homomorphism from
$\GL(M,R)$ onto
$\GL(M/M_i,R_{k-i} )$ \ $( 0\leq i\leq k)$.}
\end{quote}

Denote by  $G_i$ \ $(0\leq i\leq k)$ the normal subgroup of $G =\GL(M,R)$ which consists of all
automorphisms $v$ of $M$ that induce the identity automorphism on $M/M_i $. 
Thus $G_i$ consists of all $v\in G$ such that $v(x) -x$ lies in $M_i$
for all $x\in M$. We examine these normal subgroups, starting with
$G_1$. 

First we note that for $v\in G_1$ and for any $x=\pi y\in \pi M$ it holds that
$v(x)= x$ (because $v(x) -x =\pi (v(y )-y) =0)$. Denote by $G_1'$ the normal subgroup of $G_1$
which consists of all $v\in G_1$ for which also $v(x)= x$ for $x\in M_1$. For $v\in G_1'$ we therefore have $v(x)= x + w( x)$, where
$ w$ is an $R$-linear mapping from $M$ to $M_1$ which is zero on
$M_1+ \pi M$. One can easily see, that $G_1'$ is an Abelian group, which 
is isomorphic to  the additive group of maps $w$. Now both
$M/\pi M$ and $M_1$ can be viewed as vector spaces over the field $R_1$. These
vector spaces have the same dimension, viz.\ $a_1 + a_2+\cdots+ a_k$. An
$R$-linear map from $M$ to $M_1$ which 
is zero on $\pi M$ determines a linear mapping from the vector space $M/\pi M$
to
the vector space $M_1$. Conversely, one sees without difficulty, that to a
linear mapping from the vector space $M/\pi M$ to the vector space $M_1$ belongs an
$R$-linear mapping from $M$ to $M_1$, which is zero on $\pi M$.
This shows that $G_1'$ is isomorphic to the additive group of the
linear maps from the vector space $M/\pi M$ to the vector space $M_1$
which are zero on the subspace $M_1 + \pi M/\pi M$ of dimension $a_1$ over $R_1$.
 So 
$G_1'$ is isomorphic to the additive group of matrices with $a_2+ \cdots + a_k$ rows and
with $ a_1 + a_2 +\cdots + a_k$ columns and with entries from $R_1$. Since $R_1$, 
viewed as a vector space over $K$, has dimension $g$, it follows from the
above that $G_1'$ is the direct sum of $(a_2+ \cdots + a_k)(a_1 + a_2+ \cdots + a_k)g$ groups
isomorphic
to the additive group $K^+$ of $K$. 

We now consider $G_1 /G_1'$.
To each $v\in G_1$ belongs an automorphism $v_1$ of $M_1$ such that $v_1 (x)=
 x$ for $x\in M_1\cap \pi M$. Conversely,  given such an automorphism $v_1$ of $M_1$,
then one may associate a $v$ to it by $v(e_1^q) = v_1(e_1^q)$, $v(e_p^q)= e_p^q$ \ $(2\leq p\leq k)$.
If $v\in G_1'$, then $v_1$ is the identity automorphism of $M_1$. It follows 
from the isomorphism
 theorem that $G_1/G_1'$ is isomorphic to the group of automorphisms $v_1$ 
of $M_1$ which leave each element of $M_1\cap \pi M$ fixed. Now $M_1$ is the direct 
sum of $M_1\cap \pi M$ and a module  $V_1'$ (namely the direct sum of the $Re_1^q$).
One may view $M_1$, $M_1\cap \pi M$
and $V_1'$ as vector spaces over $R_1$. To each $v' \in G_1 /G_1'$ belongs an automorphism
of $V_1'$, since $V_1'$
is isomorphic to  $M_1 /M_1\cap \pi M$. Conversely, one can extend an automorphism of $V_1'$ directly to
an automorphism of $M_1$ that fixes every element 
of $M_1\cap \pi M$. For an automorphism $v_1'$ of $M_1$ that fixes every 
element  of $M_1\cap \pi M$  and induces the identity automorphism on $V_1'$
one has
$v_1'(x) = x + w_1'(x)$, where $w_1'$ 
is a linear map from $M_1$ to $M_1\cap \pi M$ which is zero on $M_1\cap \pi M$. So if we denote by 
$G_1''$ the normal subgroup of $G_1$, consisting of the $v\in G_1$, for which also
$v(x) -x\in M_1\cap\pi M$ for $x\in M_1$, it follows from the foregoing that $G_1/G_1''$
is isomorphic to the group of automorphisms of the vector space $V_1'$ and
that $G_1''/G_1'$ is isomorphic to the additive group of matrices with $a_1$ rows and
$a_2 + \cdots + a_k$ columns and entries from $R_1$.  Therefore $G_1/G_1''$ is
isomorphic to  $\GL_{a_1}(R_1 )$, and $G_1''/G_1'$  is direct sum of $a_1(a 2 + \cdots +
a_k)g$ groups   isomorphic to  $K^+$.

Now $G/G_i$ \ $( 1\leq i\leq k-1)$ is isomorphic to 
$\GL(M/M_i ,R_{k-i})$. The normal subgroup of $\GL(M/M_i ,R_{k-i})$ which plays the
role of $G_1$ in $G$, is  isomorphic to  $G_{i + 1}/G_i$. One thus arrives at the
following result\footnote{See also \cite{D1} for the case of an
algebraically closed field $K$.}:

\begin{quote}
\emph{
$G =\GL(M,R)$ has a series of normal subgroups $G = G_k \supset G_{k-1}\supset G_{k-1}''\supset G_{k-1}'\supset\cdots
\supset G_{1}\supset G_{1}''\supset G_1'\supset G_0 = \{1\}$, such that 
\begin{itemize}
\item $G/G_{k-1}$  is isomorphic to  $\GL_{a_k}(L)$; 
\item $G_i/G_i''$ is isomorphic to  $\GL_{a_i} (L) $; \item
$ G_i''/G_i'$  is the direct sum of $a_i ( a_{i+1} +\cdots
 + a_k) g$  groups  isomorphic to  $K^+$; \item
$G_i'/G_{i-1}$ is the direct sum  of $(a_i+\cdots + a_k)(a_{i+l} + \cdots+ a_k)g$ groups 
isomorphic to  $K^+$ \ $( 1\leq i\leq k-1)$.
\end{itemize}}
\end{quote} 

 Here \begin{itemize}
\item $G_i$ is the normal subgroup consisting of the $v\in G$ for which $v(x) -x\in M_i$ \ $(x\in M)$;
\item   $G_i''$ is the normal subgroup consisting of the $v\in G$ for which $v(x) -x\in M_i$ \ $(x\in M)$
 and $ v(x) -x\in M_{i-1}  + M_i\cap\pi M$ \ $(x\in M_i )$;
\item $G_i'$ is the normal subgroup consiting of the $v\in G$ for which $v(x) -x\in M_i$ \ $(x\in M)$
 and $ v(x) -x\in M_{i-1}  $ \ $(x\in M_i )$;
\end{itemize}

As in the case $i=1$, $M_i/M_{i-1}$ can be understood as a
vector space 
over $R_1$. This vector space is the direct sum of $(M_{i-1} + M_i\cap\pi M)/M_{i-1}$ 
and a vector space $V_i'$ of dimension $a_i$ over $R_1$, which is thus isomorphic to
$V_i = (M_i/M_{i-1} )/(M_{i- 1} + M_i\cap\pi M /M_{i-1})$ \ $(1\leq i\leq k)$. 

\section{Investigation of the classification of conjugacy classes in the symplectic group, reduction
to two cases.} \label{3} In the sequel we suppose that the
characteristic of $K$ is not equal to $2$. 

We assume to be given\footnote{See for the notions that follow here p.~3--5 in the book \cite{D} by J.~Dieudonn\'e.} 
 a skew-symmetric form $(x,
y)$ on the
$n$-dimensional vector space $E$ over $K$
i.e., a map from $E\times E$ to $K$ for which $(x_1 + x_2 ,y) = (x_1
,y) + 
(x_2 ,y)$,  \ $(\lambda x,y) = \lambda(x,y)$ \ $(\lambda \in K)$,  and $(x,y) = -(y,x)$. It follows that also 
$(x,y_1 + y_2 ) = (x,y_1) + (x,y_2 )$,  \ $(x,\lambda y) = \lambda (x,y)$. 
We furthermore suppose that this form is non-degenerate, i.e., that from $( x,
y) = 0$ for all $y\in E$ it follows that $x = 0$. The dimension of
$E$ must then be even. 

A linear transformation $u$ of $E$ for which it holds that 
\begin{equation}\label{symplectic}(u(x),u(y)) =(x,y)\quad (x,\ y\in E)
\end{equation}
is called a symplectic transformation. One sees immediately that $u$ is an
invertible transformation (from $u(x) = 0$ it follows that $(x,y) = 0$ for all $y\in E$,
so that $x = 0$). The symplectic transformations form a subgroup
$\Sp_n(K)$ of $\GL_n(K)$. We want to investigate the classification of conjugacy classes in $\Sp_n (K)$.

If $u$ is a symplectic transformation of $E$, then it follows from (\ref{symplectic}) that 
$(u(x),y) = (x,u^{- 1}(y))$ for all $x$ and $y$ from $E$. One finds that
\begin{equation}\label{pol}(f(u)(x),y) = (x,f(u^{-1} )(y))\quad (x,\ y\in E,\ f\in K[X]). 
\end{equation}
 In the sequel, for $f$ a polynomial of degree $g$, we denote by $\bar f$ the polynomial which
is determined by $ \bar f(X) = X^gf(X^{-1})$. One then has $\bar{\bar f}=  f$. Let $d$ denote the
minimal polynomial of the symplectic transformation $u$. According to (\ref{pol}) one then has
$(d(u)(x),y) = (x,d(u^{-1})(y)) =0$. From this
it follows that for all $y\in E$, \ $d(u^{-1})(y) = 0$, so also $\bar d(u)(y)=
 0$. Therefore, $\bar d$ must be a multiple of $d$ and since the degree of $\bar d$
is not greater than that of $d$, we see that $\bar d$ is of the form $ \alpha d$ \ $(\alpha\in K)$. Because $d=\overline{ (\bar d)} =
 \alpha \bar d = \alpha^2d$ we have $\alpha=\pm  1$. Write $d$ as a product of powers  of
mutually distinct irreducible polynomials  in $K[X]$, say $d=\prod_{i=1}^r  p_i^{k_i}$.
Then $\prod_{i=1}^r  p_i^{k_i}=\pm\prod_{i=1}^r  \bar p_i^{k_i}$. Thus, the polynomial $\bar p_i$ occurs
amongst 
the scalar multiples\myfootnote{\emph{``scalar multiples'' added; the leading coefficient of $\bar p_i$ is $p_i(0)$.}} 
of the polynomials $p_j$ \ $( 1\leq j \leq r)$. It follows  that the decomposition of $d$ 
can be written as $d = \gamma \prod_{i=1}^s ( p_i\bar p_i)^{k_i}\prod_{i=1}^t  q_i^{m_i}$, where the $p_i$ \ $(1\leq i\leq s)$  
 are mutually distinct irreducible polynomials for which $p_i\neq \pm \bar p_i$ 
and where the $q_i$ \ $(1\leq i\leq t)$ are mutually distinct irreducible
polynomials for which $q_i =\pm \bar q_i$, and where $\gamma \in K$. 

Denote by  $E_i$ the subspace of $E$ generated by the $x\in E$ 
for which either $p_i^{k_i}(u)(x) = 0$ or $\bar p_i^{k_i} (u)(x) = 0$ \ $(1\leq i\leq s)$. Denote by  $F_j$ the
subspace of $E$ consisting of the $x\in E$ for which it holds that $q_j^{m_j}(u)(x)= 0$ \
$(1\leq j\leq t)$. The subspaces $E_i$ and $F_j$ are invariant under $u$, and $E$ is the direct sum of 
the $E_i$ and $F_j$ (see section~\ref{1}). Furthermore,  $E_i$ 
is  the direct sum of two subspaces $E_i^1$ and $E_i^2$, consisting of 
the $x\in E$ for which $p_i^{k_i} (u) (x) = 0$, respectively $\bar p_i^{k_i} (u) (x) = 0$. The subspaces $ E_i^1$ and $E_i^2$
are also transformed into themselves by $u$. 
Suppose that the vectors $x$
and $y$ of $E$ lie in two different subspaces amongst the $E_i$ and $F_j$. Then
there are polynomials $f$ and $g$ such that $f(u)(x) = 0$, $g(u)(y)= 0$ and such that
$\bar f$ and $g$ are relatively prime ($f$ and $g$ are certain polynomials
$p_i^{k_i}\bar p_i^{k_i}$
or $q_j^{m_j})$.\myfootnote{\emph{We removed a comma between $p_i^{k_i}$ and $\bar p_i^{k_i}$.}}
Then there are two polynomials $h$ and $n$ such that $h\bar f + ng = 1$.
According to (\ref{pol}), we have $(x,y) = (x,(h\bar f+ng)(u)(y)) = (x,h\bar f(u)(y)) = (f(u)(x),u^a
h(u)(y))=  0$ \ ($a$ is the degree of $f$). We recall\footnote{\cite[p.~5]{D}} that a 
subspace $G$ of $E$ is called isotropic if there is an $x\neq0$ in $G$, such that for all
$y\in G$ it holds that $(x,y)= 0$. The subspaces $E_i$ and $F_j$ that we found here
are not isotropic. For instance, if for any $x\in E_1$ we have $ (x,y) =0$ 
for all $y\in E_1 $, then, since $E$ is direct sum of the $E_i$ and $F_j$ and since
$(x,z)=0$ for for $z\in E_i$ \ $(i\neq1)$ and for $z\in F_j$, we get $(x,y) = 0$ for all $y\in E$,
from which it follows that $x =0$. We have thus found 

\begin{quote}
\emph{$E$ is the direct sum of the
non-isotropic subspaces $E_i$ and $F_j$ \ $( 1\leq i\leq s,\  1\leq j\leq t)$.}
\end{quote}

  We want to
investigate when two symplectic transformations $u$ and $u'$ lie in
the same conjugacy class of $\Sp_n(K)$, i.e., when there is a symplectic
transformation $w$ such that $u' =wuw^{-1}$. According to section~\ref{1}, this requires
that $u$ and $u'$ have the same minimal polynomial and that the subspaces $E_i$, $F_j$ 
belonging to $u$ have the same dimension
as the subspaces $E_i'$, $F_j'$ belonging to $u'$ \ $(1\leq i\leq s,\ 1\leq j\leq t) $.
 From a theorem of Dieudonn\'e\footnote{\cite[proposition 2, p.~6]{D}} it follows that there is then a
symplectic transformation $v$ which maps the subspace $E_i'$ onto $E_i$ and
the subspace $F_j'$ onto $F_j$ \ $(1\leq i\leq s,\ 1\leq j\leq t)$. Then $vu' v^{-1}$ is  a
symplectic transformation which lies in the same conjugacy class of $\Sp_n (K)$ as
$u'$ and involves the same subspaces $E_i$ and $F_j$ as $u$. We 
assume in the sequel --which thus may be assumed without objection in the investigation of the classification of conjugacy classes in
$\Sp_n(K)$-- that $u'$ already has this
property. 

Now assume that there is a symplectic transformation $w$
so
that $u'=wuw^{-1}$. Then for an $x\in F_j$ we have $q_j^{m_j}(u)(x)= 0$, so that $0 = 
q_j^{m_j}(w^{-1}u'w)(x)=w^{-1}q_j^{m_j}(u')(w(x))$. It follows that $w(x)\in F_j$. So $w$
transforms $F_j$ \ $( 1\leq j\leq t)$ into  itself, and in an analogous way we
prove that $w$ also transforms $E_i$ \ $(1\leq i\leq s)$ into itself. 
As the restriction of $w$ to a subspace $E_i$ or $F_j$ is a symplectic 
transformation belonging to the restriction of the given form $(x, y)$
to that subspace, we see that if $u$ and $u'$ are conjugate in the
symplectic group of $E$, then the restrictions of $u$ and $u'$  to the 
subspaces $E_i$, \ $F_j$ \ $( 1\leq i\leq s,\ 1\leq j\leq t)$ are conjugate in the symplectic
groups of those subspaces. The converse is immediately apparent. 

The
investigation of whether or not two symplectic
transformations $u$ and $u'$ in the symplectic group are conjugate  is thus
reduced to the case where $u$ and $u'$  both have minimal polynomials\myfootnote{\emph{One no longer insists that a minimal polynomial must have leading coefficient $1$.}} $(p_i\bar p_i)^{k_i}$
 or $q_j^{m_j}$. We shall therefore in the 
sequel assume that $E$ itself is one of the subspaces $E_i$ or $F_j$.

\section{First case, reduction to the one dealt with in section~\ref1.}\label{4}

We start with the simplest case, namely that $u$ and $u'$ are symplectic
transformations of $E$ with minimal 
polynomial 
$(p\bar p)^k$,  where $p$ is an
irreducible polynomial of degree $g$ such that $p\neq\pm\bar p$. Recall that
$E$ is direct sum of subspaces $E^1$ and $E^2$, consisting of the $x\in E$ 
for which $p^k(u)(x) = 0$, respectively $\bar p^k(u)(x) = 0$. 
The subspaces are both transformed by $u$ into
themselves. One can find two polynomials $h$ and $f$ such that
$h\bar p^k+fp^k=1$. It follows that for $x,y\in E^1$ we have $(x,y) = 
(x,(h\bar p^k+fp^k)(u)(y)) = (x,h\bar p^k(u),y)= (p^k(u)(x),u^{gk}h(u)(y)) =0$.
Similarly, $(x,y) = 0$ for $x,y\in E^2$. Thus $ E^1$ and $E^2$ are totally 
isotropic subspaces of $E$,\footnote{\cite[p.~6]{D}} i.e., subspaces of $E$ such that the
restriction of $(x,y)$ to those subspaces is zero. Suppose that $n =  2t$ is the
dimension of $E$. A totally isotropic subspace has 
a dimension $\leq t$.\footnote{\cite[p.~7]{D}\label{twice}} Since $E$ is the direct sum of $E^1$ and $E^2$, one has 
$\dim(E^1 ) + \dim(E^2 ) = n$. Thus $\dim(E^1)=  \dim(E^2) = t$, so that $E^1$ and $E^2$
are both maximal totally isotropic subspaces. There is then  a
basis ($e_i$) of $E$ such that for
$ x =\sum_{i=1}^n\xi_ie_i$, $y =\sum_{i=1}^n\eta_ie_i$ one has 
$(x,y)=\sum_{i=1}^t(\xi_i\eta_{t+i}-\eta_i\xi_{t+i})$ (a basis with this property is called a 
symplectic basis) and further such that $e_1, \dots , e_t$ forms a basis of $E^1$ and $e_{t+1}, \dots , e_{2t}$
forms a basis of $E^2$.  With respect to the
basis $(e_i)$ the map $u$ is represented  by a matrix of the form 
\begin{equation}\label{hyp}
\begin{pmatrix}A&0\cr 0& \check A\end{pmatrix}
\end{equation}
where $A$ is invertible and $\check A$ is  the contragredient matrix (the inverse of the
transpose matrix).
${}^{\ref{twice}}$ 
From this one finds 
that the restriction of $u$ to $E^2$ is determined by the restriction to $E^1$ (and
vice versa). 

Similarly, to $u'$  belong two subspaces $E_1^1$ and $E_1^2$ for which the same holds.
So there is a symplectic basis $(e_i')$ of $E$ such that  $e_1', \dots , e_t'$ forms a basis of $E_1^1$ and $e_{t+1}', \dots , e_{2t}'$
forms a basis of $E_1^2$. The map $e_i \mapsto e_i'$ defines a symplectic transformation $w$. Then
$w^{-1}u'w$ is a symplectic transformation that lies in the same conjugacy class of the
symplectic group as $u'$ and that uses the same $E^1$, $E^2$
subspaces as belong to $u$. We  assume henceforth that $u'$ itself already had this
property: we suppose that $E_1^1 = E^1$, $E_1^2 = E^2$, and that $e_i=
e_i'$.

 With respect to the basis $(e_i)$, the map $u'$ is represented by a
matrix of the form $\begin{pmatrix}A_1&0\cr 0& \check A_1\end{pmatrix}$.
 If $u$ and $u'$ are conjugate 
in $\Sp_n (K)$, then the restrictions of $u$ and $u'$  to $E^1$  $(E^2)$ are
conjugate in the group $\GL(E^1,K)$  ($GL(E^2,K)$): from $u' = tut^{- 1}$ it follows
 that $p^k(u')(t(x)) = tp^k(u)(x) = 0$ for $x\in E^1$, so that also $t(x)\in E^1$. 
Conversely, if the restrictions of $u$ and $u'$  to $E^1$ are conjugate in $\GL(E^1,K)$,
then
there is an invertible matrix $B$ such that $A_1 = BAB^{- 1}$. The matrix 
$\begin{pmatrix}B&0\cr 0& \check B\end{pmatrix}$
then gives a symplectic transformation $t$, for which $u'= tut^{-1}$. 
It follows  that the conjugacy class of the symplectic transformation $u$ of
$E$ is completely determined by the conjugacy class  in
$\GL(E^1,K)$ of the restriction of $u$ to $E^1$.  

Finally, given a symplectic basis $(e_i)$ of $E$, let us
call $E^1$ \ $(E^2)$ the maximal totally isotropic subspaces of $E$ with
basis vectors $e_1,\dots, e_t$ \ $(e_{t+ 1},\dots ,e_{2t})$. Then starting from
a linear transformation of $E^1$ with minimal polynomial $p^k$, one can find a
symplectic transformation of $E$ with minimal polynomial $(p\bar p)^k$ which
on $E^1$ induces this linear transformation. Thus, given a conjugacy class of linear transformations 
of $E^1$ with
minimal polynomial $p^k$,  there are
symplectic transformations $u$ of $E$ such that the restriction from $u$
to $E^1$ lies in the class. From section~\ref{1} it follows that: 

\begin{quote}
\emph{To a symplectic
transformation $u$ of $E$ with minimal polynomial $(p\bar p)^k$, where $p$ is an
is an irreducible polynomial of degree $g$ such that $p\neq\pm p$, belong a
number of integers $a_i\geq0$ \ $ (1\leq i \leq k,\ a_k>0)$, 
such that $g\sum_{i=1}^ki a_i=n/2$; the conjugacy class of $u$ in the symplectic group is
unambiguously determined by the numbers $a_i$.}
\end{quote}

If $u$ is a symplectic transformation of $E$ with minimal polynomial
$(p\bar p)^k$ and if $v$ is a symplectic transformation of $E$ commuting with $u$, then, 
since the minimal polynomials 
of the restrictions of $u$ to the subspaces $E^1$, $E^2$ are relatively prime, 
 $v$ must transform these subspaces into themselves according to section 
\ref{2}. Then
$v$ is represented by a matrix of the form (\ref{hyp}) relative to the basis $(e_i)$. 
It follows  that the group of symplectic transformations commuting with $u$ 
is isomorphic to the group of invertible linear transformations of $E^1$ commuting with the restriction of $u$ 
to $E^1$. 

\section{Second case,
preparation.}\label{5} We are now going to investigate the conjugacy of two symplectic
transformations in the second case mentioned in section~\ref3, that of a
minimal polynomial $q^m$ with $q = \pm\bar q$. For this some preparations are
required. 

We assume that in the $n$-dimensional vector space $E$ over the
field $K$ with characteristic $\neq2$ the symplectic transformation $u$ has a
minimal polynomial $q^m$, where $q$ is an irreducible polynomial of
degree $h$ such that $q = \alpha \bar q$ with $\alpha^2=1$. Suppose $h\geq2$. Then $q( 1)\neq0$, $q(-1)\neq0$.
 From $q = \alpha \bar q$ it follows that $q(1)=\alpha q(1)$, \ $q(-1) =\alpha(-1)^h   q(-1)$, which implies
(as ${\rm char}(K)\neq2$) that $\alpha = 1$ and $h$ is even. So either $h$ is even and $q =\bar q$, or $h=
1$. In the latter case, $q(X) =X \pm 1$. 

For the sequel it is convenient 
to know that the 
residue class ring $K[X]/(q^m)$ is isomorphic to the residue class ring with respect to the
principal ideal $(q^m)$ in the ring $K[X]_q$ of rational fractions $\frac{u}{v}$ \ ($u,\ v\in K[X]$
with $v \not\equiv0 \bmod {(q)}$). (This is a well-known fact; incidentally, it is also
easy to verify.) Henceforth, by $R$ we mean the
residue class ring $K[X]_q /(q^m)$. 

We take a closer look at $R$. The assignment $\frac{u(X)}{v(X)}\mapsto\frac{u(X^{-1})}{v(X^{-1})}$ 
  gives an 
automorphism of $K [X]_q$ because from $v \not\equiv0\bmod ( q)$ it follows, using $q =\pm\bar q$, that
$v(\frac1X)\not\equiv 0\bmod (q)$. The ideal $(q^m)$ is mapped onto itself. So this induces an automorphism
$\rho\mapsto \bar \rho$ \ $(\rho\in R)$ of $R$.\myfootnote{\emph{not to be confused with $\bar f$ defined  for $f\in K[X]$.}}
If $\kappa$  is the canonical homomorphism from $K[X]_q$ onto $R$, then for $\rho=\kappa(\frac{u(X)}{v(X)})$  one gets
$\bar\rho=\kappa(\frac{u(X^{-1})}{v(X^{-1})})$.
It
is clear that $\rho\mapsto\bar\rho$ is an involutory automorphism, i.e., that $\overline{(\bar\rho)} = \rho$. 
Let us again put $\xi=\kappa(X)$. Then $\bar\xi=\frac1\xi$. 
 Since $R$ is generated by $\xi$,
it follows that $\rho\mapsto\bar\rho$ is the identity automorphism only if
 $\xi=\frac1\xi$, so if $X^2-1\equiv0\bmod{(q^m)}$. This can only be the case if $q^m= X\pm 1$
(as the characteristic of $K$ is different from 2). 

 Denote by  $S$ the subring of $R$ consisting of the $\sigma\in R$
for which $\bar\sigma=\sigma$. Except in the cases $q^m = X \pm1$, we have $S \neq R$.
Every $\sigma\in S$ can be written in the form $\rho + \bar\rho$ with $\rho\in R$ (e.g., with $\rho=\sigma/2$).
Conversely, an element $\rho + \bar\rho$ of $R$ lies in $S$. In particular, $\eta=\xi+\frac1\xi$ lies in $S$ and $\xi$
 satisfies the equation $\xi^2-\eta\xi+1=0$  with coefficients from $S$. When $R$ is a field, $R$ has  rank $2$
over $S$. 

Furthermore, we will need linear forms $l$ on the
vector space $R$ over $K$, i.e.\ $K$-linear maps from $R$ to $K$. These $l$
form a vector space $R'$ over $K$ (the dual vector space).\footnote{For the  notions of linear form and dual space, see \cite[$\mathsection$ 4]{B1}.
} To
any linear form $l$ and any $\lambda\in R$ one can associate a linear
form $l_\lambda$, which is determined by $l_\lambda(\rho) = l(\lambda\rho) $ \ $(\rho\in R)$. 
The map $\lambda\mapsto l_\lambda$ is a
linear map from $R$ into $R'$. The $\lambda$ for which $l_\lambda = 0$, form an
ideal $\mathfrak A(l)$ in $R$. We call $l$ degenerate if $\mathfrak A(l)$ does not consist of the
zero element of $R$ only. If $l$ is non-degenerate, then $\lambda\mapsto l_\lambda$
defines an injective linear mapping from $R$ \emph{onto} $R'$. Furthermore, 
we note that for degenerate $l$ the ideal $\mathfrak A(l)$ contains the minimal ideal
$(\pi^{m-1})$ of $R$, where $\pi = \kappa(q)$. It follows immediately that $l$ is a degenerate linear form if and only if $l(\rho) = 0$ for $\rho \equiv0\bmod{(
\pi^{m-1})}$. 

To each $l$ belongs a linear form $\bar l$, 
determined by $\bar l(\rho) = l(\bar\rho)$ \ $(\rho\in R)$. When $l$ is non-degenerate, $\bar l=l_\epsilon$  holds
with an $\epsilon\in R$. From $l(\bar\rho) = l(\epsilon\rho)$ it follows that $l(\rho) = l(\epsilon\bar\rho) = l(\epsilon\bar\epsilon\rho)$. 
So $\epsilon\bar\epsilon = 1$. We will show that a non-degenerate form $l$
can be found for which $\epsilon = \pm1$. 

If there is an element $\rho\neq0$ from
the ideal $(\pi^{m-1})$ which does not
lie in $S$, then $S$ is a proper subspace of the vector space $R$. One can thus
find a linear form $l$ such that $l(\rho) \neq 0$, $l(\sigma) = 0$ for $\sigma\in S$. Then
$l(\tau + \bar\tau) = 0$ for all $\tau\in R$. This $l$ is non-degenerate and the corresponding
$\epsilon$ is $-1$. 

Now suppose, that $\rho = \bar \rho$ for all $\rho$ from
$(\pi^{m-1})$.  Then for
all polynomials $f\in K [X]$ $$f(X)(q(X))^{m-1} \equiv f(\frac1X)(q(\frac1X))^{m- 1} \pmod {q^m},$$ from which
follows\myfootnote{\emph{Recall that $q(X)=\alpha X^{h}q(\frac1X)$, with $\alpha=\pm1$. In the subsequent computation we have inserted
$\alpha^{m-1}$ at appropriate places.}} $$\alpha^{m-1}X^{h(m-1)}f(X) \equiv f(\frac1X) \pmod q.$$
So $$\alpha^{m-1}X^{h(m-1)}\equiv 1 \pmod q, \quad[\mbox{take }f =1],$$  and $$\alpha^{m-1}X^{h(m-1)} X\equiv \frac1X \pmod q, \quad[\mbox{take } f = X].$$

However, from these congruences it follows that $$\alpha^{m-1}X^{h(m-1)}(X^2 -1) \equiv0 \pmod q,$$
which is only possible if $q = X\pm 1$. Thus, if $h\geq2$ there is a non-degenerate $l$, such that $l ( \rho+\bar\rho) = 0$ \ $ ( \rho\in R)$. 

If $h = 1$, then $q = X \pm1$.
The ideal $(\pi^{m-1})$ now consists of all multiples $\delta \pi^{m-1}$ with $\delta\in K$.
Furthermore, $\pi+\bar\pi=\xi\pm1+\frac1\xi\pm1\equiv 0 \bmod {(\pi^2)}$.
So for even $m>0$ we have $\bar\pi^{m-1}=-\pi^{m-1}$.
Since now $\pi^{m-1}$ does not in lie in $S$, one can again find an $l$ such that $l(\rho +\bar\rho)=0$ \ 
$(\rho\in R)$, \ $l(\pi^{m-1})\neq0$.
If, on the other hand, $m$ is odd, then this is not possible,
because now $\bar\pi^{m-1}=\pi^{m-1}$, so that from $l(\rho +\bar\rho)=0$ for all $\rho\in R$ it follows that
$l(\pi^{m-1})= 0$. However, now the set of $\rho-\bar\rho$ \ $ (\rho\in R)$ is a proper
subspace of the vector space $R$, and it does not contain $\pi^{m-1}$. One can therefore find a
linear form $l$
for which it holds true that $l(\rho -\bar\rho) = 0$, \ $l(\pi^{m-1} ) \neq0$. This $l$ is non-degenerate
and the corresponding $\epsilon$ equals $+1$. Thus, we have proved 

\begin{quote}
\emph{There exist non-degenerate linear forms $l$ on $R$ for which it holds true that $l(\bar\rho)= l(\epsilon\rho)$ \ $(\rho\in R)$, with
$\epsilon=-1$ if $h\geq2$ and with $\epsilon= -(-1)^m$ if $h= 1$.}
\end{quote} 

\section{Second case, reduction to
the investigation of the equivalence of Hermitian forms on a module
over a ring.}\label{6} We use the notations of section~\ref{5}. As in section~\ref{1}
one can give $E$ the structure of an $R$-module. Call this module  $N$. The skew-symmetric form $( x, y)$ given on
$E$  defines a map from $N\times 
N$ to $ K$, which we also denote by $(x, y)$. This map is $K$-linear
in $x$  (resp.\ $y$) for fixed $y$ (resp.\ $x$), and $(x,y) = -(y,x)$. Furthermore, according to
formula (\ref{pol}) of section~\ref{3}, $$(\rho x,y)=(x,\bar\rho y) \quad (x,\ y\in N,\ \rho\in R).$$ Suppose that $M$ is an
$R$-module such that there exists a bijective $R$-linear map $\phi$
from $M$ onto $N$. Like $N$ (see section~\ref{1}), $M$ is a direct sum of modules $Re_i^j$ \ $( 1\leq i\leq m,\ 1\leq j\leq b_j
)$. 
If $\phi_1$ is a
second $R$-linear mapping from $M$ onto $N$, then we have $\phi_1 = \phi t$, where $t$
is an automorphism of $M$. 

For $x,\ y\in M$ we have  the $K$-linear form $( \rho\phi(x) ,\phi(y))$ 
 on $R$. If $l$ is a fixed non-degenerate form on $R$
such that $\bar l = l_\epsilon$, then one can write 
\begin{equation}\label{l}
(\rho\phi(x),\phi(y)) = l(\rho f(x,y))\quad\mbox{ with }
f(x,y)\in R.
\end{equation}  Then $l(\rho f(x_1 +x_2 ,y)) = (\rho\phi(x_1+x_2 ),\phi(y)) = (\rho\phi(x_1
),\phi(y)) +  (\rho\phi(x_2
),\phi(y))
 =l (\rho f(x_1 ,y)) + l(\rho f(x_2 ,y))$, from which it follows that 
$l(\rho(f(x_1 +x_2 ,y)-f(x_1 ,y)-f(x_2,y))) = 0$ for all $\rho\in R$. As $l$ is not 
degenerate it follows that $f(x_1 +x_2 ,y) = f(x_1 ,y) + f(x_2 ,y)$ \ $ (x_1,\ x_2\in M)$.

Also, $l(\rho \sigma f(x,y)) = (\rho\sigma\phi(x),\phi (y)) =( \rho \phi (\sigma x),\phi (y)) = l(\rho f(\sigma x,y))$,
so that $f(\sigma x,y) = \sigma f(x,y)$. Finally, $l(\rho f(x,y)) = 
(\rho \phi (x),\phi (y)) = -(\phi (y),\rho \phi (x)) = -(\bar\rho\phi(y),\phi (x)) = -l(\bar\rho f(y,x)) = 
-l(\epsilon\rho \overline{f(y,x)})$, so that $f(x,y)= -\epsilon\overline{f(y,x)}$. We see that $f$ satisfies
\begin{equation} \label{f}\begin{cases}f(x_1+x_2
,y) = f(x_1 ,y) + f(x_2 ,y),\\
f(\rho x,y) = \rho f(x,y),\\ f(x,y) = -\epsilon \overline{f(y,x)}.
\end{cases}
\end{equation}

 From these relations it follows
easily that $f(x,y_1+y_2 ) = f(x,y_1) + f(x,y_2)$ and $f(x,\rho y) =\bar\rho f(x,y)$. In section~\ref5 we saw that we can take $l$ such that
 $\epsilon = \pm1$. By
analogy with the case where $R$ is a field, we will call an $f$ with the properties (\ref{f}) a symmetric or
skew-symmetric Hermitian form on the $R$-module $M$, 
depending on whether $-\epsilon = +1$ or $-\epsilon= -1$. We   call such a form degenerate, if there is an $x\neq0$
for which $f(x, y) = 0$ for all $y\in M$. The form $f$ found here is not
degenerate; this follows directly from (\ref{l}) since the given
skew-symmetric form on $E$ is non-degenerate. 

We call two forms $f$ and $g$ satisfying (\ref{f})
 (with the same $\epsilon$) equivalent if there is an automorphism $t$
of $M$ for which $f(t(x),t(y))= g(x,y)$ holds. Since the map $\phi$  from $M$
onto $N$ is unambiguously determined by $u$ except for an automorphism of $M$, the
found form $f$ on $M$ is determined up to equivalence. 

Let $M$ be a given $R$-module (with the same dimension over $K$ as $E$)
and let there be given an equivalence class of non-degenerate forms $f$ on $M$.
We will show that
there exists a symplectic transformations $u$ of $E$ with minimal polynomial $q^m$
and such that $u$ determines on $M$ the equivalence class of  $f$.
So suppose $\tau$ is a bijective $K$-linear
mapping from $M$ onto $E$, and suppose $f(x,y)$ is a non-degenerate form on $M$
satisfying (\ref{f}). Then we get
a skew-symmetric form $k$ on $E$ determined by  $k(\tau(x), \tau(y)) = l(f(x,y))$;
 it is clear that $k$ is a
bilinear form, and that $k$ is a skew-symmetric form follows from
$k(\tau(x),\tau(y)) = l(f(x,y)) = -l(\epsilon \overline{f(y,x)})= -l(f(y,x))  =
-k(\tau(y),\tau(x))$. Since $f$ on $M$ is non-degenerate, $k$ on $E$ is non-degenerate. Now
two non-degenerate skew-symmetric forms on $E$ are equivalent.\footnote{For any skew-symmetric form $s$
there is a basis with respect to which $s$ is given by 
$s(x,y)=\sum_{i=1}^t(\xi_i\eta_{t+i}-\xi_{t+i}\eta_i)$ \ $(n=2t)$. See \cite[p.~5]{D}.} There
is thus an invertible $K$-linear transformation $t$ of $E$ such that $k(x,y) =
(t(x),t(y))$ \ $(x,\ y\in E)$. Setting $u'(\tau(x))  = \tau(\xi x)$ determines a linear
transformation $u'$ of $E$ for which $k(u'\tau(x),u'\tau(y))  = l(f(\xi x,\xi y))=
l(f(x,y)) = k(\tau(x),\tau(y))$. 
Thus, $u'$ is  a symplectic transformation belonging to the form $k$. 
Then $u = tu't^{-1}$ is a symplectic transformation belonging to the given form $(x,y)$  on $E$.
 The  ring of linear
transformations of $E$ generated by $u$ is isomorphic to $R$, and can be considered as a ring
of operators on $E$. Thus one can give $E$ the structure of an $R$-module.
Call that $R$-module $N$.  Then $\phi = t\tau$ is  a mapping from $M$ onto $N$ for which $(\phi(x),\phi(y)) = l(f(x,y))$. It is an $R$-linear map: for
$x\in M$ we have $\phi(\xi x) = t\tau(\xi x)  = tu'\tau(x) = ut\tau(x) = u\phi (x) = \xi\phi(x)$;
it is clear that $\phi(x+y)
= \phi(x) + \phi(y)$. \\Thus, $u$ is a symplectic transformation
of $E$ that determines on $M$ the equivalence class of $f$.

\

We
investigate when two symplectic transformations $u$ and $u'$  of $E$
with minimal polynomial $q^m$ lie in the same conjugacy class of the symplectic group. 
Starting from $u$ or from $u'$  one can define on $E$ the structure of an
$R$-module. Call the modules one 
 gets $N$ and $N'$.  According to section~\ref{1}, $u$ and $u'$  can only lie in the same
class of the symplectic group, if $N$ and $N'$ are isomorphic. We will therefore assume this. There are then bijective
$R$-linear maps $\phi$, $\phi'$ from the fixed $R$-module $M$ to $N$,  respectively $N'$, 
 such that  $u(\phi (x)) =
\phi (\xi x)$, \ $u' (\phi'  (x)) = \phi' (\xi x) $ \ $(x\in M)$. Furthermore, there are forms $f$ and $f'$ on
$M$ which satisfy (\ref{f}) and for which $(\phi (x),\phi (y)) =  l(f(x,y))$ 
and  $(\phi' (x),\phi' (y)) = l(f' (x,y))$ \ $(x,y\in M)$. Now suppose that there is a
symplectic transformation $w$ of $E$ such that $u' = w^{- 1}uw$. Then 
$t=\phi^{-1}w\phi'$ is an $R$-linear mapping from $M$ onto itself: $t(x+y) =  t(x) +
t(y)$ is obvious and furthermore $t(\xi x) = \phi^{ -1}w\phi' (\xi x) = 
\phi^{ -1}wu' \phi'(x)=\phi^{ -1}uw\phi' (x)=\xi \phi^{-1}w \phi' (x)=\xi t(x) $ \ $(x\in M)$, from which it 
follows that $t(\rho x) = \rho t(x) $ \ $(\rho\in R,\ x\in M)$. That $w$ is a symplectic transformation
is expressed by $(w\phi' (x),w\phi' (y)) = (\phi' (x),\phi' (y))$ viz.\ by
$(\phi t(x),\phi t(y)) = (\phi' (x),\phi' (y))$ \ $(x,\ y\in M)$. It follows that $l( f( t(x),
t(y))) = l( f' (x, y))$. Because $l( \rho f( t(x), t(y)) )= l( \rho f' (x, y))$
for all $\rho\in R$, we have $f(t(x),t(y)) = f' (x,y)$. Thus $f$ and $f'$ are equivalent on $M$
 if $u$ and $u'$  lie in the same class of the symplectic group.
Conversely, if there is an automorphism $t$ of $M$ such that $f(t(x),t(y)) = f'
(x,y)$ then it is easy to verify that $w 
= \phi t(\phi')^{-1}$ is a symplectic transformation of $E$ for which $u' =
w^{-1}uw$.

We summarize what we found in this section: 

\begin{quote}
\emph{To a symplectic
transformation $u$ of $E$ with minimal polynomial $q^m$, where $q$ is an irreducible
polynomial such that $q = \pm \bar q$, belongs a module $M$  over $R= K[X]_q /(q^m)$
with the same dimension over $K$ as $E$ and an equivalence class of  non-degenerate symmetric or 
skew-symmetric Hermitian forms on $M$.  For a given module $M$ and a given class of non-degenerate forms
on $M$ there is at least one $u$. Two symplectic transformations $u$ and $u'$ with the same $M$
are conjugate in the symplectic group if and only if the same class of Hermitian forms on $M$ belongs
to $u$ and $u'$.}
\end{quote} 

\section{Investigation of the equivalence of two Hermitian forms on a module
over a ring.}\label{7}
 The result of section~\ref{6} requires to look further into
the equivalence of two symmetric or skew-symmetric
Hermitian forms on the $R$-module $M$. We use the same notations as in the previous section.
The ring $R$ has an automorphism $\rho\mapsto\bar\rho$. Since this
automorphism maps the ideal $(\pi)$ onto itself, there is an induced automorphism of $R_i=R/(\pi^i)$. 
 We denote this also by
$\rho\mapsto\bar\rho$ \ $( \rho\in R_i)$. 

We begin with the case where $q$ has degree $h\geq2$. Then $h$ is
 even (see section~\ref{5}). For $\pi_1 = \xi^{-h/2}\pi$ one has $\pi_1=\bar\pi_1$, and $\pi_1$ also generates 
  the ideal $(\pi)$. Since in formula (\ref{f}) of section~\ref6 we can now take $\epsilon= -1$, we need only consider 
  the case of a symmetric
Hermitian form on $M$ for our investigation of symplectic
transformations.  But first we prove:

\begin{quote}
\emph{If $f$ is a non-degenerate symmetric or skew-symmetric Hermitian
form on $M$, then for every $R$-linear mapping $x\mapsto \lambda(x)$ from $M$ into
$R$ there is $y\in M$ such that $\lambda(x)=f(x,y)$ for all $x\in M$.}
\end{quote}

 Recall that $M$ is
direct sum of modules $Re_i^j$ \ $(1\leq i\leq m,\ 1\leq j\leq b_i )$. The $R$-module $M'$ of 
linear forms $\lambda(x)$ (the dual module) is then the direct sum of the modules $R(e_i^j)'$ \
$( 1\leq i\leq m,\ 1\leq j\leq b_i)$, where the linear form $( e_i^j)'$
is defined by $(e_i^j)'(e_k^h)=0$  if $i\neq k$ or $j\neq h$,  \ $(e_i^j)'(e_i^j) = \pi^{m-i}$.
This shows that $M'$ is isomorphic to $M$. On the other hand, 
one also obtains a mapping from $M$ onto a 
 submodule  of $M'$  that is isomorphic to  $M$ by assigning to each $y\in M$ the linear form
$f(x,y)$. But $M'$ cannot
be isomorphic to  a proper submodule, since a proper subspace of $M'$ as a vector space over $K$ has a
dimension smaller than the dimension of $M'$. It follows  that the
submodule  of the $f(x,y)$ is identical with $M'$. 

\

We now show that a
Hermitian form $f$ on $M$ also induces Hermitian forms on the residue class modules that were introduced in section~\ref{2}. 
 We use the notations of section~\ref2. Let $f$ be a
symmetric Hermitian form. 

\begin{quote}
\emph{(a) On the $R_{m-i}$-module $M/M_i$ \ $(1\leq i\leq m-1)$, $f$ determines  a
symmetric Hermitian form $f_i$, which is non-degenerate if $f$ is non-degenerate.}
\end{quote}

Denote by  $\phi_{ m-i}$ and $\psi_i$ the canonical homomorphisms of $R$ onto $R_{m-i}$ and of $M$
onto $M/M_i$ respectively. Then define $f_i$ by $f_i (\psi_i (x), \psi_i  (y)) = \phi_{m-i} (f(x,
y))$ \ $(x, y\in M)$.  One checks without difficulty that $f_i$ is a symmetric Hermitian form 
on the $R_{m-i}$-module $M/M_i$. If $f_i(\psi_i (x),\psi_i  (y)) = 0$ for an 
$x\in M$ and for all $y\in M$, then one has $f(x,y) \equiv 0 \pmod {\pi^{m-i}}$ for all $y\in M$, 
from which it follows that $f(\pi^i x,y) = 0$ for all $y\in M$, so that $\pi^i x = 0$ and  $\psi_i  (x) =
0$ if $f$ is non-degenerate on $M$. 

\begin{quote}
\emph{(b) On the vector space $V_i = (M_i /M_{i-1}
 )/(M_{i-1} +M_i\cap\pi M/M_{i-1} )$ \ $ (1\leq i\leq m)$ over $R_1$,  $f$ determines a symmetric
Hermitian form $f_i'$, which is non-degenerate if $f$ is non-degenerate.}
\end{quote} 

 We
    first note that for $x\in M_i$ one has $f(x,y)\equiv  0 \pmod{\pi^{m-i}} $. 
Denote by  $\chi_{m-i}$ the $R$-linear mapping from the ideal $(\pi^{m-i})$ onto $R_i$, 
determined by $\chi_{m-i}(\pi_1^{m-i}\rho) =\phi_i  (\rho)$ \ $(\rho\in R)$. Denote by  $\psi_i'$ the
canonical homomorphism from $M_i$ onto $V_i$ ($\psi_i'$ is the map obtained by
first mapping $M_i$ canonically onto $M_i/M_{i - 1}$ and then $M/M_{i - 1}$  onto $V_i$). 
For $\rho\in R$ and $x\in M_i$, one then gets $\psi_i'(\rho x) =\phi_1 (\rho)\psi_i (x)$ and $\chi_{m-i} (\overline{\pi_1^{m-i}\rho}) 
=\chi_{m-i} (\pi_1^{m-i}\bar\rho) =\phi_1(\bar\rho) = \overline{\phi_1(\rho)}  = \overline{\chi_{m-i}(\pi_1^{m-i}\rho)}$. Now define $f_i'(\psi_i'(x),\psi_i'(y))=\chi_{m-i}(f(x,y))$  for $x,\ y
\in M_i$. Since  $\chi_{m-i}(f(x,y))= 0$
when
$x$ or $y$ lies in $M_{i-1} + M_i\cap\pi M$, this determines a mapping from $V_i\times V_i$ into
$R_1$. One effortlessly verifies that $f_i'$ is a symmetric
Hermitian form on $V_i$.

We now show that $f_i'$ is non-degenerate if $f$ is non-degenerate. We
assume that $f$ is non-degenerate. First we check for which $x\in M$
one has $f(y,x) = 0$ for all $y\in M_i$.  As  $\pi_1^{m-i}y\in M_i$ for all 
$y\in M$, we get  for such an $x$ that $\pi_1^{m-i} f(y,x) = 0$, from which it follows that $f(y,x) \equiv 0 \bmod{(\pi^i)}$.
Say $f(y, x) = \pi_1^i g(y)$. Then $g(y)$ is determined modulo 
$\pi^{m-i}$, and for $y\in M_i$ we have $g(y) \equiv 0 \bmod{(\pi^{m-i})}$.  Put $g_1(\psi_i (y)) = \phi_{m-i}(g(y))$. Then $g_1$  is  a linear
form on the $R_{m-i}$-module $M/M_i$.  Since, according to (a), $f_i$ is non-degenerate on
$M/M_i$ one can find  $z\in M$ such that $f_i(\psi_i (y),\psi_i(z)) = g_1(\psi_i
(y))$ \ $(y\in M)$, viz.\ $\phi_{m-i}(f(y,z))= \phi_{m-i} (g(y))$, from which it  follows that 
$f(y,\pi_1^iz)=f(y,x)$ for all $y\in M$. This is only possible if $x = \pi_1^iz$, i.e., if
$x\in \pi^iM$. If, conversely, $x\in \pi^iM$, then it is clear that $f(y,x) = 0$ for
all $y\in M_i$. 

Now suppose that for certain $x\in M_i$ and for all $y\in M_i$ one has
$ f_i'(\psi_i'(y),\psi_i'(x))= 0$. Then $f(y,x)\equiv0 \bmod{(\pi^{m-i+1})}$ and $ f(y,\pi^{i- 1}x) =0$
 for all $y\in M_i$. It follows  that $\pi^{i-1}x\in \pi^i M$. Since $x\in M_i$, also $x\in M_{i-1}+M_i\cap \pi M$,
so that $\psi_i' (x) = 0$. Therefore $f_i'$ is  non-degenerate on $V_i$ if $f$
is non-degenerate.

\begin{quote}
\emph{(c) If all the forms $f_i'$ \ $(1\leq i\leq m)$ on $V_i$ are non-degenerate, then $f$ 
is non-degenerate on $M$.}
\end{quote} 

 For $m = 1$ this assertion is true (then $M= V_1$).
Suppose it has already been proved for modules $M$ over a ring $K[X]_q/(q^n)$ with
$n<m$. Just like for $M$ one has the vector spaces $V_i$ and the forms $f_i'$ \ $ ( 1\leq i\leq m)$.
And one has for $ M/M_a$ \ $( 1\leq i\leq m-a)$  vector spaces $W_i$ over   $R_1$  and forms $f_i''$ on
these vector spaces $( 1\leq i\leq m-a)$. The vector space $W_i$ has the same dimension over $R_1$ as
$V_{i+a}$ and from the definition of the forms it further follows easily that there
is a linear map $t_i$ from $W_i$ onto $V_{i +a}$ such
that $f_i''( x,y) =f_{i+a}' (t_i (x),t_i (y))$ \ $(x,y\in W_i )$. 

Now suppose that $f(x,y) = 0$ for some $x\in M$ and for all $y\in M$. Then $f_1 (\psi_1
(x),\psi_1(y)) = 0$ for all $\psi_1 (y)\in M/M_1$. Since, according to the 
induction hypothesis, $f_1$ is non-degenerate on $M/M_1 $, one has 
$\psi_1(x) = 0$, i.e., $x\in M_1$.  From $f(x,y) = 0$ for all $y\in M_1$ it further follows
that
$f_1' (\psi_1'
(x),\psi_1'(y)) = 0$  for all $\psi_1'(y)\in V_1$. Since $f_1'$ is non-degenerate, one has
$\psi_1'(x) = 0$, i.e. $x\in M_1\cap \pi M$. Say $x=\pi z$. Then $f(z,y)\equiv0\bmod{(\pi^{m-1})}$
for all $y\in M$. So $f_1(\psi_1(z),\psi_1(y)) = 0$ for all 
$\psi_1(y)\in M/M_1 $. It follows that $\psi_ 1(z) = 0$, or $z\in M_1 $. But then
$x =\pi z= 0$. Therefore $f$ is  non-degenerate on $M$. 

\

Now suppose that $f$ and $g$ are two non-degenerate symmetric forms on $M$. 
To $f$ and $g$ belong non-degenerate symmetric forms $ f_i$ and $g_i$ on $M/M_i$, 
as well as $f_i'$ and $g_i'$ on $V_i$. We
prove then:

\begin{quote}
\emph{For the forms $f$ and $g$ to be equivalent on $M$ it is necessary and sufficient that the forms $f_i'$ and $g_i'$ are equivalent on $V_i$ \ $(1\leq i\leq m)$.}
\end{quote}

 The necessity can be seen
immediately  by observing that an automorphism of $M$ which transforms $f$ into $g$ induces on $V_i$ 
an automorphism that transforms $f_i'$ into $g_i'$.

We prove the other part of the assertion by complete
induction on $m$. For $m = 1$, we have $M = V_1$ and the assertion is true. Assume
that for modules $M$ over the ring $K [X]_q / ( q^n)$ with $n<m$ the
assertion has already been proved. On the $R_{m -1}$-module $M/M_1$ the forms $f_1$, $g_1$ are then 
 equivalent. Since $M_1$ is direct sum  of $M_1\cap\pi M$ and a vector space $V_1'$
isomorphic to  $V_1$  (see section~\ref2), and $f(x,y) = 0$ for $x\in \pi M$, $y\in M_1$,
the equivalence of $f$
and $g$ on $M_1$
 follows from the equivalence of $f_1'$ and $g_1'$.
 Since, according to the result of section~\ref2, there is an automorphism
 of $M$ which induces given  automorphisms on $M/M_1$ and on $V_1$, there exists
an automorphism $u_1$ of $M$ such, that 
\begin{equation}\label{fgrad}
\begin{cases}f(u_1(x),u_1 (y)) \equiv g(x,y)\pmod { \pi^{m-1}}\quad
(x,\ y\in M),\\
f(u_1 (x),u_1 (y)) = g(x,y)\quad (x,\ y\in M_1).
\end{cases}
\end{equation}
 In order to derive the claim from (\ref{fgrad}), we need three auxiliary results.
 
 \

I) If $h$ and $k$ are two non-degenerate symmetric Hermitian forms on
$M$ such that, for certain $j$ \ $(1\leq j\leq m-1)$, we have\myfootnote{\emph{We will repeatedly use \emph{or}, where the original used  \emph{and} \textup{[}in Dutch: \emph{en}\textup{]}.}}
\begin{equation*}
\begin{cases}h(x,y) \equiv k(x,y) \pmod { \pi^{m-1}}\quad
(x,\ y\in M),\\
h(x,y) = k(x,y)\quad (x,\ y\in M_j \mbox{ or }x\in M_{j-1},\ y\in M),
\end{cases}
\end{equation*}
 then there is an
automorphism $v$ of $M$ such that
\begin{equation}\label{hgrad}
\begin{cases}h(v(x),v(y)) \equiv k(x,y) \pmod { \pi^{m-1}}\quad
(x,\ y\in M),\\
h(v(x),v(y))= k(x,y)\quad (x\in M_j,\   y\in M).
\end{cases}
\end{equation}

We try to find $v$ with $v(x) = x + \pi_1^{j-1}w(x)$, where $w$
is an $R$-linear mapping from $M$ into $M_j$ such that $\pi_1^{j-1}w(x)  = 0$ for $x\in M_j$.
It is easy to see that such $v$ is an automorphism of $M$ and that the
first relation in (\ref{hgrad}) is satisfied. It follows from the other relation in (\ref{hgrad}), that $w$
must satisfy $$\pi_1^{j-1}h(x,w(y)) =k(x,y) -h(x,y) \quad(x\in M_j,\ y\in M).$$
Write $k(x,y) - h(x,y) = \pi_1^{m-1}p(x,y)$. Then it follows  that for $x\in M_j$
it must hold true that  
\begin{equation} \label{hprime}h_j'(\psi_j'(x),\psi_j'(w(y))) = \phi_1(p(x,y)).
\end{equation}
For $x\in M_{j-1} + M_j\cap\pi M$, we have $k(x,y) -h(x,y) =0$, so $p(x,y) \equiv 0 \bmod {(\pi)}$. This
shows that  $\lambda(\psi_j'(x)) = \phi_1 (p(x,y))$ defines a linear form $\lambda$ on $V_j$. This form  is zero if $y\in M_j$. 
Recall the elements $e_p^q$ \ 
$(1\leq p\leq m,\  1\leq q\leq b_p)$  of $M$ such that $M$ is the direct sum
of the
modules $Re_p^q$. Since $h_j'$ on $V_j$ is not
degenerate, we can  define for any element $e_p^q$ an element
$w(e_p^q)$  so that (\ref{hprime}) is satisfied. One can further suppose that 
$w(e_p^q) = 0$ for $1\leq p\leq j$. The $w(e_p^q)$ determine an endomorphism $w$ of $M$. 
One easily verifies that $\pi_1^{j-1}w(x) = 0$ for $x\in M_j$.
Thus, for the $v$ derived from this $w$, (\ref{hgrad}) is satisfied. 

\

II) If $h$ and $k$ are two non-degenerate symmetric Hermitian forms on
$M$ such that for certain $j$ \ $(1\leq j\leq m-1)$ one has
\begin{equation*}
\begin{cases}h(x,y) \equiv k(x,y)
\pmod{\pi^{m-1}}  \quad (x,\ y\in M),\\
h(x,y)= k(x,y)\quad (x\in M_j,\   y\in M),
\end{cases}
\end{equation*}
then there is an
automorphism $v'$ of $M$ such that 
\begin{equation}\label{vgrad}
\begin{cases}h(v'(x),v'(y)) \equiv k(x,y) \pmod { \pi^{m-1}}\quad
(x,\ y\in M),\\
h(v'(x),v'(y)) = k(x,y) \quad (x,\ y\in M_{j+1}\mbox{ or } x\in M_{j},\ y\in M).
\end{cases}
\end{equation}

We try to find $v'$ with $v' (x) = x + \pi_1^jw'(x)$, where $w'$ is a
linear mapping from $M$ into $M_{j+1}$ such that $w' (e_p^q) = 0$  for $p\neq j
+1$. One easily sees that such $v'$ is an automorphism of $M$ and that
then the first relation in (\ref{vgrad}) is already satisfied. The other relation in (\ref{vgrad})
gives
\begin{equation}\label{w'}
\begin{aligned}
 \pi_1^jh(&w'(x),y) + h(x,w'(y)) + \pi_1^{2j}h(w'(x),w'(y)) =\\ &k(x,y)
-h(x,y) \quad (x,y\in M_{j+1} \mbox{ or } x\in M_j,\ y\in M).
\end{aligned}
\end{equation}

Since $w' (x)\in M_{j+1} $, we have $h(w' (x),w' (y))\equiv 0\pmod { \pi^{m-j-1}}$. It follows that
the third term on the left hand side of (\ref{w'}) is zero. Furthermore, it is immediately apparent
that (\ref{w'}) is satisfied for $x\in M_j$, $y\in M$. Now put $k(x,y) 
-h(x,y) = \pi_1^{m-1}p'(x,y)$. Then it follows from (\ref{w'}) that $w'$ needs to satisfy
\begin{equation}\label{p'}
\begin{aligned}
 h_{j+1}'(\psi_{j+1}'(&w'(x)),\psi_{j+1}'(y)) + h_{j+1}'(\psi_{j+1}'(x),\psi_{j+1}'(w'(y)) )  =\\ &
\phi_1(p'(x,y)) \quad (x,\ y\in M_{j+1}).
\end{aligned}
\end{equation}
If one puts $q(\psi_{j+1}'(x),\psi_{j+1}'(y))=\phi_1(p'(x,y))$, then $q$ is a
symmetric Hermitian form on $V_{j+i}$. Auxiliary result 
II) then follows from 

\

III) If $V$ is a
vector space over the field $R_1$,   $h$ is a non-degenerate symmetric
Hermitian form on $V$ and $n$ a symmetric Hermitian form on $V$, then there is
a linear transformation $t$ of $V$ such that $$h(t(x),y) + h(x,t(y)) =
n(x,y).$$

This is easy to prove: for any fixed $y\in V$, we see that $\frac12n(x,y)$ is a
linear form on $V$. One can then find a linear transformation $t$ with the
intended property  from $\frac12n(x,y)=  h(x,t(y))$. 

\

The
assertion about the equivalence of the forms $f$ and $g$ that we wanted to prove 
is a direct consequence of 
I) and II): starting from (\ref{fgrad}) we can, by
application of 
I) and II), find an automorphism $u_i$  of
$M$ for each $i$ \ $ (1\leq i\leq m)$  such that $f(u_i (x),u_i (y)) = g(x,y)$ is satisfied for
$x,\ y\in M_i$ and for $x\in M_{i- 1},\ y\in M$. Then $u_m$ is such 
such that $f(u_m (x),u_m (y)) = g(x,y)$ for all $x$ and $y$ from  $M$.

\

Finally, we prove: 

\begin{quote}
\emph{Given forms $f_i'$ on the $V_i$ \ $(1\leq i\leq m)$, there is a symmetric Hermitian form $f$ that induces
them.
 }
\end{quote}

Namely, suppose that
$ x=\sum_{p=1}^m\sum_{q=1}^{b_p}\xi_p^qe_p^q$ and $y=  \sum_{p=1}^m\sum_{q=1}^{b_p}\eta_p^qe_p^q$
are two elements 
 of $M$. Then for a symmetric Hermitian form $f$, one has $f(x, y) = 
\sum_{p,\ r=1}^m\sum_{q=1}^{b_p}\sum_{s=1}^{b_r}\alpha_{pr}^{qs}\xi_p^q\bar\eta_r^s
$ with $\alpha_{pr}^{qs}=\bar\alpha_{rp}^{sq} $, \ $\pi^p\alpha_{pr}^{qs}=\pi^r\alpha_{pr}^{qs} =0$.
One  easily sees 
that if $\alpha_{ii}^{qs}=\pi_1^{m-i}\beta_{ii}^{qs}$, one has for $x,\ y\in M$ that 
$f_i'(\psi_i'(x),\psi_i'(y))=  \phi_1(\sum_{q,\ s=1}^{b_i}\beta_{ii}^{qs}\xi_i^q\bar\eta_i^s)$.
 By appropriately choosing the  $\alpha_{ii}^{qs}$, it is possible to achieve that $f_i'$ becomes equal to a given form.

\section{Continuation of section~\ref7.}\label{8} In section~\ref7, we examined the equivalence of
the forms occurring in the study of symplectic transformations
in the case where the polynomial $q$ has degree $h\geq2$. 
We now consider the case $h =1$.
In the study of symplectic
transformations with a minimal polynomial $q^m$ we now get to deal with a form $f$ satisfying
relations (\ref{f}) of section~\ref6, where $-\epsilon= (-1)^m$. (See the result of
section~\ref5). Furthermore, we now have $\pi+\bar\pi\equiv 0 \bmod{(\pi^2)}$, and the automorphism that is induced on $R_1$
is the identity automorphism (as $\xi-\bar\xi\equiv0\bmod{ (\pi)}$, one has
$ \rho-\bar\rho \equiv0\bmod{ (\pi)}$ for each $\rho\in R$). 

\

Our  starting point  is a Hermitian form $f$
on $M$ which  satisfies (\ref{f}) of section~\ref6 with $ \epsilon= \pm1$ (we
thus do not yet suppose that $- \epsilon=  (-1)^m$). As in section~\ref7 
one can then define  a form $f_i$ on $M/M_i$ which is symmetric
or skew-symmetric precisely if $f$ is so, and which is non-degenerate if $f$ is non-degenerate. 

The
forms $f_i'$ on $V_i$ \ $(1\leq i\leq m)$, however, must now be defined differently.
Denote by  $\chi_{m-i}'$ the $R$-linear mapping from the ideal $(\pi^{m-i})$ into $R_1$ 
determined by $\chi_{m-i}'(\pi^{m-i}\rho)=\phi_1(\rho)$ \ $ (\rho\in R)$. Then, since $\pi+\bar\pi\equiv 0\bmod{(\pi^2)}$, we get
$\chi_{m-i}'\overline{(\pi^{m-i}\rho)}=\chi_{m-i}'((-\pi)^{m-i}\bar\rho+\pi^{m-i+1}\sigma)$, 
 with certain $\sigma\in R$
for each $\rho\in R$. It follows that $\chi_{m-i}'\overline{(\pi^{m-i}\rho)}=(-1)^{m-i}\phi_1(\bar\rho)=(-1)^{m-i}\phi_1(\rho)=
(-1)^{m-i}\chi_{m-i}'(\pi^{m-i}\rho)$.
 Now define $f_i'$ by
$f_i'(\psi_i'(x),\psi_i'(y))=\chi_{m-i}'({ f(x,y)})$ \ $(x,\ y\in M_i)$.
 Then
$f_i'(\psi_i'(x),\psi_i'(y)) = \chi_{m-i}'(\overline{-\epsilon f(y,x)}) = -\epsilon(-1)^{m-i}\chi_{m-i}'(f(y,x)) =
-\epsilon(-1)^{m-i}f_i'(\psi_i'(y),\psi_i'(x))$. Thus $f_i'$ is now   a symmetric or
skew-symmetric bilinear form on $V_i$, depending on whether $-\epsilon(-1)^{m-i} $ equals  $+1$ or
$-1$. One proves in the same way as in section~\ref7 that $f_i'$ is non-degenerate
if $f$ is non-degenerate and that $f$ is non-degenerate  if all $f_i'$
are non-degenerate $(1\leq i\leq m)$. 

For the case $-\epsilon = (-1)^m$, which
we have to deal with for symplectic transformations, we note that it follows that the $f_i'$
with odd $i$ are non-degenerate skew-symmetric bilinear forms on the
vector space $V_i$ of dimension $b_i$ over $R_1$.  So $b_i$ is even if $i$ is odd.

Again, for the equivalence
of the non-degenerate symmetric or skew-symmetric Hermitian
forms $f$ and $g$ on $M$ it is necessary and sufficient that the forms $f_i'$ and $g_i'$ be equivalent on $V_i$.
This can be deduced from claims 
I) and II) of section~\ref7 (which
must then be slightly modified in obvious ways).
I) is proved in the same way for the present case. The proof of 
II)
however, becomes slightly different. Instead
of (\ref{w'}) one gets the relation 
$\pi^jh(w'(x),y) + \bar\pi^jh(x,w'(y))=
k(x,y) -h(x,y)$ \ $(x,\ y\in M_{j+1})$.
Put $k(x,y) -h(x,y) = \pi^{m-1}p'(x,y)$. Then $\pi^{m- 1}p'(x,y) = 
-\overline{\epsilon\pi^{m- 1}p'(y,x)}=-{\epsilon(-\pi)^{m- 1}p'(y,x)}$, from which it follows that $\phi_1(p'(x,y)) =
-\epsilon(-1)^{m- 1}\phi_1(p'(y,x))$.  In place of (\ref{p'}) we now need
\begin{equation}\label{p'bis}
\begin{aligned}
 h_{j+1}'(\psi_{j+1}'(&w'(x)),\psi_{j+1}'(y)) +(-1)^j h_{j+1}'(\psi_{j+1}'(x),\psi_{j+1}'(w'(y) ) ) =\\ &
\phi_1(p'(x,y)) \quad (x,\ y\in M_{j+1}).
\end{aligned}
\end{equation}
Here $h_{j +l}'$ is such that $h_{j+1}' (x, y) = -\epsilon(-1)^{m-j - 1}h_{j+1}' (y,x)$ \ $(x,\ y\in V_{j+1})$.
 For 
 III) we must now substitute 
 
 \
 
III${}'$) If $V$ is a vector space
over the field $R_1$, $h$ is a non-degenerate bilinear form on $V$ such
that $h(x,y) = \epsilon_1 h(y,x)$ \ $ (\epsilon_1^2 = 1)$,  $n$ is a bilinear form on $V$ such that
$n(x,y) =\epsilon_2n(y,x)$ \ $ (\epsilon_2^2 =1)$, then there is a linear transformation $t$ of $V$
such that $$h(t(x),y) + \epsilon_1\epsilon_2 h(x,t(y)) = n(x,y).$$ 
(Thus in (\ref{p'bis}) one has $\epsilon_1=-\epsilon
(-1)^{m-j-1}$, \ $\epsilon_2 = -\epsilon
(-1)^{m-1}$.)

 The proof is
easy: one may determine $t$ by $\frac12n(x,y) =  h(x,t(y))$. 

\

It is proved in the same
way as in section~\ref7 that there are forms $f$ which induce given
forms $f_i'$ on the $V_i$.

\section{The classification of conjugacy classes in the symplectic group.}\label{9}
 We summarize what we found in
sections~\ref3 through \ref8 about the classification of conjugacy classes in the group $\Sp_n(K)$
(where the characteristic of $K$ differs from 2) : 

\begin{quote}
\emph{Every conjugacy class
of $\Sp_n(K)$ is unambiguously determined by \\
\mbox{ \ }
(a) a number of mutually distinct irreducible polynomials $p_i \in K[X]$,
different from $X$, of degrees $g_i$ and such that $p_i \neq \pm\bar p_i$ \ $ (1\leq i\leq s)$, integers
 $k_i >0$ and integers $a_j^i\geq0$ \ $ (1\leq i\leq s,\ 1\leq j \leq k_i,\   a_{k_i}^i 
>0)$;\\ 
\mbox{ \ }
(b) a number of mutually distinct irreducible polynomials $q_i \in K[X]$
of degrees $h_i\geq2$ and such that $q_i=\bar q_i$ \ $(1\leq i\leq t)$, integers $m_i>0$, 
integers $b_j^i\geq0$ \ $ (1\leq i\leq t,\ 1\leq j\leq m_i, \ b_{m_i}^i>0)$ and equivalence classes $K_j^i$ of 
non-degenerate symmetric Hermitian forms on the vector space of
dimension $b_j^i$ over the field $L_i=K[X]/(q_i)$ \ $ (1\leq i\leq t,\ 1\leq j\leq m_i)$;\\
\mbox{ \ \ }(c) powers $(X+1)^{m_+}$ and $(X-1)^{m_-}$ of the polynomials $X+1$ and $X-1$ respectively,  integers
$ b_j^+\geq 0$ \ $ (1\leq j\leq m_+,\  b_{m_+}^+ > 0)$ and $b_j^- \geq0$ \ $ ( 1\leq j\leq m_- ,\  b_m^- > 0)$ such
that the $ b_j^+$ and the $b_j^-$ with odd index $j$ are even, and equivalence classes $K_j^+$ and $K_k^-$
\ $(1\leq j\leq m_+,\  1\leq k\leq m_-,$ $j$ and $k$ even$)$ of  non-degenerate
quadratic forms on the vector space of dimension $b_j^+$ and $b_k^-$ respectively over $K$;\\
where the following equation is satisfied $$2\sum_{i=1}^sg_i\sum_{j=1}^{k_i}j a_j^i +\sum_{i=1}^th_i\sum_{j=1}^{m_i}jb_j^i+
\sum_{j=1}^{m_+}jb_j^++\sum_{j=1}^{m_-}jb_j^-=n.$$}
\end{quote} 

We make a few more remarks about this result in special
cases.

1) \underline{$K$ is an algebraically closed field.} Every irreducible
polynomial then has the degree $1$, so that there are no polynomials $q$. Furthermore
two non-degenerate quadratic forms on an $n$-dimensional
vector space over $K$ are equivalent (for every form $f$ there is a 
basis $(e_i)$ such that for $ x = \sum_{i=1}^n\xi_ie_i$,\ $ y= \sum_{i=1}^n\eta_ie_i$
one has $ f(x,y) =  \sum_{i=1}^n\xi_i\eta_i$).\footnote{See for instance \cite[p.~65]{We}} 
A conjugacy class 
of the symplectic group is thus unambiguously determined by the
polynomials $p_i$ and the numbers $k_i$,\ $a_j^i$, and by the numbers $m_+$, $b_j^+$ and 
$m_-$,\ $ b_j^-$. It follows  that now two symplectic transformations lie in
the same conjugacy class of $\Sp_n(K)$ if they lie in the same class of $\GL_n(K)$.

 2) \underline{$K$ is the field of  real numbers.}\footnote{For this case the classification of conjugacy classes 
 in $\Sp(K)$ has been investigated in a different way by J.~Williamson \cite{W}.} Two non-degenerate
quadratic forms on a vector space over $K$ are 
equivalent if and only if their indices of inertia are equal.\footnote{See  \cite[p.~130]{vdW} for the concept of index of inertia.}  The same
holds for non-degenerate Hermitian forms on a vector space over the
field of complex numbers. It follows  that to each $K_j^i$
unequivocally belongs  a number $c_j^i$ \  $(0\leq c_j^i\leq b_j^i)$ and 
to each $K_j^+$, $K_k^-$ a number $c_j^+$ \ $ (0\leq c_j^+\leq b_j^+)$ resp.\ $c_k^-$  \  $(0\leq c_k^-\leq b_k^-)$. 

3) \underline{$K$ is a
finite field with characteristic  different from 2.} If $L$ is a
finite field with an involutory automorphism $\rho\mapsto\bar \rho$, then for
every non-degenerate Hermitian form $ f$ (with respect to this automorphism) on an
$n$-dimensional vector space over $ L$ one has a basis $ (e_i')$ of the
vector space such that 
for $x=\sum_{i=1}^n\xi_ie_i'$, \  $ y \sum_{i=1}^n\eta_ie_i'$ one has $f(x,y) = \sum_{i=1}^n\alpha_i\xi_i\bar\eta_i$,\footnote{\cite[p.~64]{D}}
where $\bar\alpha_i=\alpha_i\neq0$.
 Furthermore, since every $\alpha_i \in L$ with  $\bar\alpha_i=\alpha_i$ can be written 
as $\alpha_i = \beta_i\bar \beta_i$ \ $(\beta_i\in L)$\footnote{This follows from well-known
theorems about the solvability of quadratic equations in a finite
field. See, e.g., \cite[p.~175]{B2}}, there is also a basis $(e_i)$ relative to 
which $f(x,y) = \sum_{i=1}^n\xi_i\bar\eta_i$. This shows that two non-degenerate
Hermitian forms on a vector space over $L$ are equivalent. Furthermore
two non-degenerate quadratic forms on a vector space over
a finite field with characteristic different from 2  are
equivalent only if their determinants are in the same
coset of the subgroup of squares  in the multiplicative group $K^*$ of $K$.\footnote{See \cite[p.~158]{LED}.} Using the foregoing, for a finite field
$K$, it is easy to specify the classification of conjugacy classes in $\Sp_n (K)$. For $n=4$ and
$n=6$ one then finds the results of L.~E.~Dickson \cite{LED1}, \cite{LED2}. 

\section{The
symplectic transformations which commute with a given symplectic
transformation}\label{10} 
 Finally, we investigate the group
of the symplectic transformations $v$ which commute with a given symplectic
transformation $u$. With the help of what was proved in section~\ref3
it is easily seen that we can restrict ourselves to the case
where $u$ has a minimal polynomial $(p\bar p)^k$ \ $ (p\neq\pm\bar p)$ or a minimal polynomial $q^m$ \
$(q=\pm\bar q)$. The first case has already been covered at the end of section~\ref4,
so here we suppose that $u$ has  minimal polynomial $q^m$, where $q$
is an irreducible polynomial of degree $h$ such that $q=\pm\bar q$. We use
the notations of the previous sections. 

With a linear transformation $v$ commuting with $u$
belongs an automorphism $t$ of the
module  $M$ that was introduced in section~\ref6 (see section~\ref2). This $t$ is determined by $\phi t(x) = v\phi (x)$ \  $(x\in M) $. If
$v$ is a symplectic transformation, then $(\phi (t(x)),\phi (t(y))) =
(\phi (x),\phi (y))$, from which it follows that $f(t(x),t(y)) = f(x,y)$. This shows that
the group of symplectic transformations commuting with the symplectic transformation $u$
 is isomorphic to the group $U=\U(M,R, f)$ of
automorphisms $t$ of the $R$-module $M$ for which $f(t(x),t(y)) =
f(x,y)$. We investigate this group in the same way as the group
$\GL(M,R)$ (see section~\ref2). 

 As in section~\ref2
an automorphism $t$ of $M$ for which  $f(t(x),t(y)) =
f(x,y)$ induces  an automorphism $v_i$ on the $R_{m-i}$-module $M/M_i$, for which 
$f_i (v_i(x),v_i(y))=f_i(x,y)$ \ $(x,\ y\in M/M_i)$. Conversely, assuming
a $v_i$ for which this holds, one can find an automorphism $t'$ of $M$ which on $M/M_i$
induces the automorphism $v_i$ and which is such that $t'(e_p^q) = e_p^q$ \ $(1\leq p\leq i)$ (see section~\ref2). 
For $ i=1$, 
by applying 
I) and II) from section~\ref7 (resp.\ section~\ref8), one can  
construct from this $t'$ an automorphism $t$ of $M$ which on $M/M_1$ induces the 
automorphism $v_1$ and which is such that $f(t(x),t(y)) =f( x,y)$. 
When $ i>1$, one can apply this to the module $M/M_{i-1}$ and thus  find an
automorphism $v_{i-1}$ of $M/M_{i-1}$ which on $M/M_i$ induces the automorphism
$v_i$ and which is such that $f_{i-1}( v_{i-1} (x), v_{i-1} (y)) = 
f_{i-1} (x, y)$ \  $(x,\ y\in M/M_{i-1})$.  From the foregoing it then follows without
difficulty:

\begin{quote}
 \emph{Assigning to an automorphism from the group $\U(M,R,f)$ 
the automorphism it induces on $M/M_i$ gives a homomorphism
from $\U(M,R, f)$ onto $\U(M/M_i,R_{m-i},f_i)$.}
\end{quote}

 By  $G$, $G_i'$ and $G_i''$ we shall mean the groups introduced in
section~\ref 2. Then $U_i =U\cap G_i$ is the normal subgroup of
$U$ which consists of the $t\in U$ that  induce the identity automorphism on $M/M_i$. 
Furthermore, we set $U_i' = U\cap G_i'$, \ $U_i'' = U\cap G_i''$ \ $(1\leq i\leq m-1)$. We are going to examine these
normal subgroups in more detail. 

\

First, suppose that $h\geq2$. We can
then suppose that $f$ is a symmetric Hermitian form on $M$ (see
section~\ref6). Now $U_i'$ consists of all automorphisms $t$ of $M$ for which $t(x) = x
+ w(x)$ with $w(x)\in M_1$, \ $w(x) = 0$ for $ x\in M_1$, while  also 
$f( t(x), t(y)) = f(x, y),$ so $$ f(x,w(y)) + f(w(x),y) + f(w(x),w(y)) =
0.\quad (x,\ y\in M).$$

Since $w(x) = 0$ for $x\in M_1$, we have $f(x,w(y)) = 0$ for $x\in M_1$, from which it follows that 
$w(y)\in \pi M$ (see section~\ref7). However, then $f(w(x),w(y)) =0$ \ $(x,\ y\in M)$. It 
follows that $U_1'$ is isomorphic to the additive group of linear maps
$w$ from $M$ into $M_1\cap \pi M$ which are zero on $M_1$ and for which $f(x,w(y))
+ \overline{f(y,w(x))} = 0$ \ $ (x,\ y\in M)$. Such a $w$ induces a 
linear mapping $w'$ from the vector space $V = (M/\pi M)/(M_1 +\pi M/\pi M)$ over $R_1$
into the vector space $W = M_1\cap\pi M$ over $R_1$. If we denote by $\chi$ the 
canonical homomorphism from $M$ onto $V$, then  $f$ determines a mapping
$h(x,y)$ \ $(x\in V,\ y\in W)$ from $V\times W$ into $R_1$ which is linear in $x$ and  anti-linear
in $y$, where $h(\chi(x),y)$ is defined by $h(\chi(x),y) =  \chi_{m-1}(f(x,y))$ \ $(x\in M,\
y\in W$, $\chi_{m-1}$ is the linear 
mapping from $(\pi^{m-1})$ into $R_1$ that was introduced in section~\ref7). 
From $h(\chi(x),y) = 0$ for all $x\in M$
follows $f(x,y) = 0$ for all $x\in M$, and from this follows $y = 0$. For the
map $w'$ one has
\begin{equation}\label{hchi}  h(\chi(x),w'(\chi(y))) + \overline{h(\chi(y),w'(\chi(x)))} = 0.
\end{equation}
Conversely, one can easily see that given a linear mapping from
$V$ into $W$ satisfying this relation one gets an unambiguously determined $ w$. 

Both $V$
and $W$ are  vector spaces of dimension $b_2+ \cdots + b_m$ over $R_1$. Choose
 bases in $V$ and $W$. If $X$ (resp.\ $Y$) represents the matrix 
of components of the vector $x$ (resp.\ $y$) from $V$ (resp.\ $W$) relative to the basis of $V$ (resp.\ $W$),
then we have $h(x,y) = X'H\bar Y$, where $H$ is an invertible matrix with entries
from $R_1$ (here, $X$ and $Y$ are matrices with one column, and if $A$ is a matrix with entries from $R_1$, we denote 
by $A'$  the
transposed matrix and by $\bar A$ the matrix obtained from $A$
by replacing each entry $\rho\in R_1$ of the matrix by $\bar\rho$). If to the
linear transformation $w'$ belongs  a matrix $L$, then the relation (\ref{hchi})  reads
in matrix form $X'H\bar L\bar Y  +\bar Y'\bar HLX=0$. The result is that $H\bar L+ L'\bar H' =0$. Putting
$H\bar L=T$ we get
 then $T + \bar T' = 0$, and $\bar L =  H^{- 1}T$. This shows that the group consisting of the $w'$
is isomorphic to  the additive group of the skew-symmetric Hermitian
matrices with $b_2 + \cdots +b_m$ rows and columns and with entries in $R_1$. One
easily verifies, that this group of matrices is the direct sum of $(b_2
+ \cdots + b_m)^2 \frac{h}2$  groups isomorphic to  $K^+$. The same thus holds for
$U_1'$. 

To each $t\in U_1$ belongs an automorphism $t_1$ of $M_1$
for which $t_1(x) = x$ for $x\in M_1\cap \pi M$ and $f(t_1(x),t_1 (y)) = f(x,y)$ \ $ ( x,\ y\in M_1$)\ (see section~\ref2).  
Conversely, one can also find a $t$ for a $t_1$ (this again turns out to be
easy using  
I) and II) from section~\ref7). The $t$'s for which $t_1$ is the
identity automorphism are those from $U_1'$. Thus $U_1/U_1'$ is isomorphic to 
the group of automorphisms $t_1$. And $U_1''/U_1'$ is isomorphic to  the group of
automorphisms $t_1'$ of $M_1$ that satisfy $ t_1' (x) = x $ \ $(x\in M_1\cap\pi M)$, \ $t_1' (x)
-x\in M_1\cap\pi M$ \ $(x\in M_1 )$, \ $f(t_1'(x),t_1'(y)) =  f(x,y)$ \ $(x,\ y\in M_1)$. However,
if the first two relations are satisfied, then the last relation is also satisfied: this follows directly 
from $f(x_1, y_1 ) = 0$ for $x_1 \in M_1$, $y_1 \in\pi M$.
Thus $U_1''/U_1'$ is isomorphic to  $G_1''/G_1'$ and 
in section~\ref2 we have seen that this quotient group is direct sum of $b_1(b_2
 +\cdots + b_m)h$ subgroups  isomorphic to  $K^+$. Furthermore, it is easy to see that $U_1
/U_1''$ is isomorphic to  the group $\U_{b_1}(f_1',R_1)$ of the 
automorphisms of $V_1$ which leave invariant the non-degenerate symmetric Hermitian form
$f_1'$. 

In the same way as in section~\ref2, the investigation
of the groups $U_i$ with $i>1$ can be reduced to the above case. The result is (with $L =K[X]/(q)$): 

\begin{quote}
\emph{If $h\geq2$ then the group $ U
= \U(M,R,f)$ has a series of normal subgroups \\
 $U = U_m \supset U_{m-1}\supset U_{m-1}''\supset U_{m-1}'\supset\cdots
\supset U_{1}\supset U_{1}''\supset U_1'\supset U_0 = \{1\}$, such that 
\begin{itemize}\item
$U/U_{m-1}$  is isomorphic to  $\U_{b_m}(L,f_m')$; \item
$U_i/U_i''$ is isomorphic to  $\U_{b_i} (L,f_i') $; \item
$ U_i''/U_i'$  is the direct sum of $b_i ( b_{i+1} +\cdots
 + b_m) h$  groups  isomorphic to  $K^+$; \item
$U_i'/U_{i-1}$ is the direct sum  of $(b_{i+l} + \cdots+ b_m)^2\frac{h}2$ groups 
isomorphic to  $K^+$.
\end{itemize}}
\end{quote}

We now consider the case $h=1$.
The result then becomes slightly different. Assume that $f( x, y) =-\epsilon\overline {f( y, x)}$%
 \ $ (
\epsilon =\pm1)$.\myfootnote{\emph{bar added.}} In the same way as above, we find that $U_1'$ is isomorphic to
the additive group of linear maps $w$ from $M$ into $M_1\cap\pi M$ which are zero
on $M_1$ and for which  $f(x,w(y)) -\overline{ \epsilon f(y,w(x))} = 0$ \ $(x,\ y\in M)$.
Define $V$, $W$ and $\chi$ in the same way and now define $h$ by $h(\chi(x),y) =
\chi_{m-1}' (f(x,y))$,
where in this case $\chi_{m-1}'$ is the linear map
from
$(\pi^{m-1})$ into $R_1$  introduced in section~\ref8. We get that $w$ induces a linear mapping $w'$ from the
vector space $V$ into the vector space $W$, for which $$h(\chi(x),w'(\chi(y))) +
(-1)^m\epsilon h(\chi(y),w'  (\chi(x)) )= 0.$$

Introducing bases in $V$ and $W$ shows that $U_1'$ is isomorphic to the additive
group of symmetric or skew-symmetric matrices with $b_2 + \cdots+ b_m$
rows and columns and with entries in $ R_1$ according to whether $(-1)^m\epsilon = -1$ or
$+1$. It follows without difficulty that $U_1'$ is direct sum of $\frac12(b_2 + \cdots
+b_m)(b_2 + \cdots +b_m + (-1)^{m- 1}\epsilon)$ groups isomorphic to  $K^+$. Again, 
$U_1''/U_1'$ is 
isomorphic to  $G_i''/G_1'$ and $U_1/U_1''$ is isomorphic to  $\Sp_{b_1} (R_1 )$ if
$\epsilon=(-1)^{m-1}$ and with the group $\Orth_{b_1} (R_1 ,f_1')$ of orthogonal 
transformations belonging to the quadratic form $f_1'$ on $V_1$ if $\epsilon = (-1)^m$.

Using this, one can also investigate the other $U_i$. One finds,
since $R_1$ is now isomorphic to  $K$, for the case $\epsilon=(-1)^{m-1}$, 
which we have to deal with for the case of symplectic transformations: 

\begin{quote}
\emph{
 If
$h=1$ and\,\myfootnote{\emph{bar added.}} \ $f(x,y) = (-1)^m\overline{f(y,x)}$,
then the group $U = \U(M,R,f)$ has a series of
normal subgroups 
$U = U_m \supset U_{m-1}\supset U_{m-1}''\supset U_{m-1}'\supset\cdots
\supset U_{1}\supset U_{1}''\supset U_1'\supset U_0 = \{1\}$,
such that \begin{itemize}\item
$U/U_{m-1}$ is isomorphic to  $\Orth_{b_m} (K,f_m')$ if $m$ is even and to $\Sp_{b_m} (K)$ if
$m$ is odd;\item
$ U_i/U_{i}''$ is isomorphic to  $\Orth_{b_i} (K,f_l')$ if $i$ is even and to $\Sp_{b_i} (K)$
if $i$ 
is odd;\item
 $U_i''/U_i'$ is the direct sum of $b_i ( b_{i+l} + \cdots + b_m)$  groups isomorphic to  $K^+$;\item
$U_i'/U_{i-1}$ is the direct sum of $\frac12(b_{i+1}+\cdots+b_m)(b_{i+1}+\cdots+b_m+(-1)^{i+1})$
groups isomorphic to  $K^+$.
\end{itemize}}
\end{quote}

Using these results, one can find, e.g., in the case where $K$
is a finite field of characteristic different from $2$,  the order of
the normalizer of an element of the symplectic group.
For $n=4$ these numbers were calculated by L.~E.~Dickson.

\section*{Summary}
\myfootnote{\emph{originally in French.}}
In this thesis we have studied the classification of conjugacy classes of elements of the symplectic group $\Sp_n(K)$, i.e.
the group of linear transformations with $n$ variables and coefficients in the commutative field $K$ 
leaving invariant an alternating bilinear form. We have limited ourselves to the consideration of the case 
where the characteristic of $K$ is different from $2$.

The result we have arrived at can be stated as follows: each conjugacy class of the group $\Sp_n(K)$ is 
characterized by a system of invariants. These invariants are first of all, as in the case of the general linear group 
(the group of all invertible linear transformations), irreducible polynomials and systems of non-negative integers, but secondly also equivalence 
classes of certain Hermitian forms and of certain quadratic forms.

Here are some indications of how we have treated the problem. Let us denote, for a one-variable polynomial $f$ of degree $g$ 
with  coefficients in $K$, by $\bar f$ the polynomial defined by $\bar f(X) = X^gf (\frac1X)$. It is easy to see that the minimal polynomial of 
a symplectic transformation of a vector space $E$ on $K$ satisfies $f = \pm \bar f$. It is shown that we may focus on two cases: 
\\
l) $f = (p\bar p)^k$, where
$p$ is an irreducible polynomial such that $p\neq\pm\bar p$,\myfootnote{\emph{The text mistakenly said $p=\pm \bar p$.}}\\
2) $ f = q^m$, where $q$ is an irreducible polynomial such that $q = \pm\bar q$. 

In the first case,  the class of 
$u$ in the symplectic group is determined by a conjugacy class  of a general linear group. 
One can thus use  the known  theory of canonical forms of linear transformations.

In the second case we proceed as follows. By virtue of the the relation $q = \pm\bar q$, the ring $R = K[X]/(q^m)$ has an
 involutive automorphism $\rho\mapsto\bar\rho$ \  $(\rho\in R)$. One can define on $E$ a module structure 
 with respect to $R$. Letting $M$ be a copy of this $R$-module, we show that $u$ and the alternating bilinear form  
 given on the vector space $E$, define on $M$ an equivalence class of anti-symmetric Hermitian forms, or symmetric 
 Hermitian forms,
  where by such a form is meant a map $f(x,y)$ from $M\times M$ to $R$ such that 
 $f(x_1 + x_2,y)=f(x_1,y)+f(x_2,y)$, $ f(\rho x,y) = \rho f(x,y)$ \ $(\rho\in R)$, $f(x,y) =\epsilon \overline{f(y,x)}$
  \ $(\epsilon=\pm1$; one can even assume that $\epsilon=+1$, 
except in the case where $q = X\pm1$, $m$ odd). 

If two symplectic transformations $u$ and $u'$  with 
minimal polynomial $q^m$ are conjugate in the general linear group, then we can take the same 
module  $M$ for $u$ and for $u'$. We show that $u$ and $u'$  are conjugate in the symplectic group if and only if $u$ and $u'$ give the same equivalence class
of Hermitian forms on $M$. Finally, we prove that two Hermitian forms on $M$  
are equivalent if and only if certain ordinary\myfootnote{\emph{here ``ordinary'' means they take values in a field.}} Hermitian forms (if the degree of $q$ is  $>2$) or certain quadratic forms 
(if the degree of $q$ is $1$) are equivalent. 

 We have also studied the structure 
of the normalizer of an element of the symplectic group G. One finds results analogous to those found by J.~Dieudonn\'e
 in the case of the general linear group \cite{D1}. The study of the conjugacy classes of the 
 other classical groups (i.e. the orthogonal and unitary groups) can be done in a similar way. 
 We hope to return to this on another occasion.

\appendix


\section*{Propositions\myfootnote{\emph{As was customary, the PhD thesis ends with a separate list of `Propositions' (`Stellingen') without proofs, about which the author could be questioned at the oral defense.}}}
\begin{enumerate}
\item The elements of the Lie algebra belonging to a
skew-symmetric or quadratic form over a field $K$, can
be divided into classes of elements that can
be obtained from each other  with
a symplectic resp.\ orthogonal transformation. One can investigate this division into classes by the same
method as was used in this thesis in the investigation of the
classification of conjugacy classes in the symplectic group. 

\item One can define the trace of a
linear transformation of a vector space without using
a basis of the vector space. 

\item It is possible to
prove purely algebraically that any ordinary orthogonal transformation
is a product of two-dimensional rotations and of reflections. 

\item It
is probable that an irreducible representation of the finite
group $\GL_n(\mathbb F_q)$ (where $\mathbb F_q$ is the finite field with $q$ elements), is
characterized by certain invariants which bear a lot of resemblance
with the invariants characterizing the conjugacy classes of this group.

\item The
theorem proved by E.~Hecke that every irreducible representation
by matrices of the group $\PSL_2 (\mathbb F_q)$ is equivalent to a representation
by matrices all entries of which lie in the field produced by the characters, can also be proved without using
arithmetic tools. 
\begin{quote} E.~Hecke, Math.\ Ann., Bd.\ 116 (1939),
p.\ 469--510. \end{quote} 

\item Some of the irreducible 
representations of the modular group modulo $p^\lambda$ can be calculated in an algebraic way.
 \begin{quote} H.~D.~Kloosterman, Ann.\ of Math., vol. 47 (1946), p.~317--447. \end{quote} 

\item It is desirable that in the theory of partial
differential equations attention should be given to a strict
definition of the concept of integral surface.

\item The inequality derived by H.~Weyl for the powers of the two
 types of eigenvalues of a linear transformation can also
be proved using  elementary matrix and  
differential calculus. 
\begin{quote} H.~Weyl, Proc.\ Nat.\ Ac.\ Sc.\ vol.\  35 (1949), 
p.~408--411. \end{quote} 

\item Tensor calculus is not indispensable for a clear formulation of
differential geometric properties.
 
\item The
calculation by K.~Husimi and I.~Sy\^ ozi of the state sum belonging to
a planar hexagonal lattice can be 
brought into a simpler form. 
\begin{quote} K.~Husimi and  I.~Sy\^ ozi, Progr.\ Theor.\ Phys., vol.\ V (1950) p.177--186,
 I.Sy\^ ozi, ibid., p.~341-351. \end{quote} 
\end{enumerate}


\bibliographystyle{elsarticle-num} 


\end{document}